\documentclass[final,reqno]{siamltex}
\usepackage[table,usenames,dvipsnames]{xcolor}
\usepackage{tikz,pgfplots}
\usetikzlibrary{arrows}
\usepackage{placeins}
\usetikzlibrary{decorations.markings}
\usepackage{placeins}
\usepackage{multirow}
\usepackage{amsmath,amssymb}
\usepackage[square,sort&compress,comma,numbers]{natbib}
\usepackage{multirow}
\usepackage{enumerate}
\usepackage{graphicx}
\usepackage{subfigure}
\usepackage[linesnumbered,ruled,vlined]{algorithm2e}
\usepackage{floatrow}
\usepackage[left=1in, right=1in, top=1in, bottom=1in]{geometry}
\newfloatcommand{capbtabbox}{table}[][\FBwidth]
\definecolor{fullrank}{rgb}{0.8,0.2,0}
\definecolor{lowrank}{rgb}{0,0.95,0.95}
\definecolor{darkgreen}{rgb}{0,0.4,0}
\newtheorem{remark}{\emph{\textbf{Remark}}}
\newtheorem{bm}{\emph{\textbf{Benchmark}}}

\title{The Inverse Fast Multipole Method}
\subtitle{An $\mathcal{O}(N)$ fast direct solver.}
\author{Sivaram Ambikasaran\thanks{Courant Institute of Mathematical Sciences, New York University, New York, NY 10012} \and Eric Darve\thanks{Institute for Computational and Mathematical Engineering, Huang Engineering Center, Via Ortega, Stanford University, Stanford, CA 94305-4042}
}

\begin{document}
\maketitle
\pagestyle{myheadings}
\markboth{Sivaram Ambikasaran, Eric F. Darve}{The Inverse Fast Multipole Method}
\addtocounter{MaxMatrixCols}{20}
\begin{abstract}
This article introduces a new fast direct solver for linear systems arising out of wide range of applications, integral equations, multivariate statistics, radial basis interpolation, etc., to name a few. \emph{The highlight of this new fast direct solver is that the solver scales linearly in the number of unknowns in all dimensions.} The solver, termed as Inverse Fast Multipole Method (abbreviated as IFMM), works on the same data-structure as the Fast Multipole Method (abbreviated as FMM). More generally, the solver can be immediately extended to the class of hierarchical matrices, denoted as $\mathcal{H}^2$ matrices with strong admissibility criteria (weak low-rank structure), i.e., \emph{the interaction between neighboring cluster of particles is full-rank whereas the interaction between particles corresponding to well-separated clusters can be efficiently represented as a low-rank matrix}. The algorithm departs from existing approaches in the fact that throughout the algorithm the interaction corresponding to neighboring clusters are always treated as full-rank interactions. Our approach relies on two major ideas: (i) The $N \times N$ matrix arising out of FMM (from now on termed as FMM matrix) can be represented as an extended sparser matrix of size $M \times M$, where $M \approx 3N$. (ii) While solving the larger extended sparser matrix, \emph{the fill-in's that arise in the matrix blocks corresponding to well-separated clusters are hierarchically compressed}. The ordering of the equations and the unknowns in the extended sparser matrix is strongly related to the local and multipole coefficients in the FMM~\cite{greengard1987fast} and \emph{the order of elimination is different from the usual nested dissection approach}. Numerical benchmarks on $2$D manifold confirm the linear scaling of the algorithm.
\end{abstract}


\section{Introduction}
\label{section_introduction}

Large dense matrices arising out of many applications such as integral equations, interpolation, inverse problems, etc., can be efficiently represented as hierarchical matrices, which are data sparse representations of certain class of dense matrices. Detailed discussions on different hierarchical matrix structures can be found in~\cite{hackbusch1999sparse,hackbusch2000sparse, hackbusch2000sparse1, grasedyck2003construction, hackbusch2001introduction, borm2003introduction, hackbusch2002data, hackbusch2000h2, hackbusch2002h2,borm2010efficient,borm2006matrix, chandrasekaran2006fast,chandrasekaran2006fast1}. The FMM enables computing fast matrix vector products for a sub-class of these hierarchical matrices, i.e., obtaining $b = Ax$, where the matrix $A$ is a FMM matrix at a computational cost of $\mathcal{O}(N)$, given a tolerance $\epsilon$. In this article, we provide an algorithm to answer the opposite question:

``Can we construct a direct solver for the linear system $Ax=b$, given the right hand side $b$ and a specified error tolerance $\epsilon$, where $A$ is a FMM matrix at a computational cost of $\mathcal{O}(N)$?"

Traditionally, solving linear systems of the form $Ax=b$, where $A$ is a FMM matrix has been done using iterative solver. However, in recent times, there has been an increasing focus on constructing fast direct solvers for these hierarchical matrices~\cite{ambikasaran2013fast,ho2012fast,martinsson2005fast,greengard2009fast,martinsson2009fast,chandrasekaran2006fast,chandrasekaran2006fast1,aminfar2014fast,kong2011adaptive}. This endeavor has been fruitful in the case of hierarchical matrices with weak-admissibility criteria (strong low-rank structure), i.e., the interaction between non-overlapping regions is represented as a low-rank interaction. This includes the HODLR~\cite{ambikasaran2013fast,kong2011adaptive, aminfar2014fast} and HSS matrices~\cite{chandrasekaran2006fast, martinsson2005fast,ho2012fast,chandrasekaran2006fast1}. However, the weak admissibility criteria (strong low-rank structure) is restrictive. The rank of the non-overlapping clusters is more or less constant for problems arising out of $1$D manifolds but is no longer valid for singular Green's function of elliptic PDE's on $2$D and $3$D manifolds, i.e., the rank of interaction between neighboring clusters is no longer independent of $N$. For instance, the rank grows as $\mathcal{O}(\sqrt{N})$ on $2$D manifolds and grows as $\mathcal{O}(N^{2/3})$ on $3$D manifolds. Hence, the scaling of existing fast direct solvers is no longer linear in higher dimensions. This article addresses this issue by constructing a fast direct solver, which relies only on compressing the interaction corresponding to non-neighboring clusters. More precisely, the article presents a new fast direct solver for solving linear systems of the form $K \sigma = \phi$, where $K \in \mathbb{R}^{N \times N}$ is an FMM matrix, in $\mathcal{O}(N)$ operations with a controllable error $\epsilon$. The algorithm also extends to the more general class of matrices known as $\mathcal{H}^2$ matrices (hierarchical matrices with nested basis) with strong admissibility criteria (weak low-rank structure).

Below are the class of problems that the solver proposed in this article is capable of solving.
\begin{enumerate}
\item
\textbf{Integral equations:}
Integral equations arising from elliptic partial differential equations are of the form:
\begin{equation}
a(\vec{x}) \sigma(\vec{x}) + \int_{\Omega} b(\vec{x}) G\left( \Vert \vec{x} - \vec{y} \Vert \right) c(\vec{y}) \sigma(\vec{y}) dy = \phi(\vec{x})
\label{eqn_integral_equation}
\end{equation}
where $\vec{x} \in \mathbb{R}^n$, $\sigma(\vec{x}) \in \mathbb{R}^n \mapsto \mathbb{R}$ is the unknown function and $a(\vec{x}), b(\vec{x}), c(\vec{x}),\phi(\vec{x}): \mathbb{R}^{n} \mapsto \mathbb{R}$ are known functions. The function $G\left( r \right)$ is the Green's function of the underlying elliptic partial differential equation. Table~\ref{table_Greens_function} lists some examples of different Green's functions.
{\renewcommand{\arraystretch}{2}
\begin{table}[!htbp]
\begin{center}
\caption{Examples of Green's function $G(r)$}
\rowcolors{1}{}{gray!40}
\begin{tabular}{|c|c|}
\hline
\textbf{Elliptic PDE} & \textbf{Green's function}\\
\hline
Laplace/Poisson in $2$D & $\log\left(r\right)$\\
Laplace/Poisson in $3$D & $1/r$\\
Helmholtz in $2$D & $H_0^{(1)}(k r)$\\
Helmholtz in $3$D & $\dfrac{\exp(-ikr)}{r}$\\
\hline
\end{tabular}
\label{table_Greens_function}
\end{center}
\end{table} 
}

Discretizing~\eqref{eqn_integral_equation} leads us to solving a system of the form:
$$K {\sigma} = {\phi}$$
where $\phi \in \mathbb{R}^N$ is known, $\sigma \in \mathbb{R}^N$ is unknown and the matrix $K \in \mathbb{R}^{N \times N}$ possesses a FMM structure~\cite{lai2014fast}.

\begin{remark}
In case of Helmholtz's, and in general any oscillatory Green's function, the solver is applicable for only moderately high frequency problems.
\end{remark}

\item
\textbf{Dense covariance matrices:}
Dense covariance matrices~\cite{ambikasaran2013fastbayes,saibaba2012application} arising in many applications in inverse problems~\cite{ambikasaran2013large}, Kalman filtering~\cite{li2014kalman,lee2013hydrogeophysical}, statistics, machine learning, Gaussian process~\cite{,ambikasaran2014fastdet,ambikasaran2014fastsym}, etc., can be efficiently represented as FMM matrices~\cite{ambikasaran2013large}. The entries of the covariance matrix, $K \in \mathbb{R}^{N \times N}$, arise from the covariance functions, $C(r)$, i.e., $K(i,j) = C(\Vert \vec{x}_i- \vec{x}_j \Vert)$. Table~\ref{table_covariance_functions} lists some examples of different choices of the covariance functions, $C(r)$.
\begin{table}[!htbp]
\begin{center}
\caption{Examples of covariance functions $C(r)$}
\rowcolors{1}{}{gray!40}
\begin{tabular}{|c|c|}
\hline
Ornstein Uhlenbeck/ Exponential & $\exp \left(- r/a \right)$\\
Squared exponential/ Gaussian & $\exp \left(-r^2/a \right)$\\
Rational quadratic & $\dfrac1{\left( 1 + r^2 \right)^{\alpha}}$\\
Mat\'{e}rn & $\dfrac{\left(\sqrt{2 \nu} r/\rho \right)^{\nu}}{\Gamma(\nu) 2^{\nu-1}} K_{\nu} \left(\sqrt{2 \nu} r/\rho \right)$\\
\hline
\end{tabular}
\label{table_covariance_functions}
\end{center}
\end{table} 

\item
\textbf{Radial basis function interpolation:}
Radial basis function interpolation for multivariate approximation is one of the most frequently applied techniques in approximation theory to represent scattered data in multiple dimensions. The literature on interpolation using radial basis function is vast and we refer the reader to a few~\cite{wu1993local, schaback1995creating, billings2002interpolation, wang2002point, buhmann2003radial, de2007mesh}. Table~\ref{table_rbf} lists the commonly used radial basis functions and the corresponding matrices, i.e., $K(i,j) = R(\Vert \vec{x}_i - \vec{x}_j\Vert)$ can be well-represented as FMM matrices~\cite{ambikasaran2013fast,beatson1999fast,gumerov2007fast,carr2001reconstruction}.

\begin{remark}
In addition to the radial basis functions in Table~\ref{table_rbf}, the covariance functions in Table~\ref{table_covariance_functions} are also used in radial basis function interpolation.
\end{remark}
\begin{table}[!htbp]
\begin{center}
\caption{List of radial basis functions $R(r)$}
\rowcolors{1}{}{gray!40}
\begin{tabular}{|c|c|}
\hline
Multi-quadric & $(1+r^2/a^2)^{1/2}$\\
Inverse multi-quadric & $(1+r^2/a^2)^{-1/2}$\\
Poly-harmonic spline I & $r^{2k+1}$\\
Poly-harmonic spline II & $r^{2k} \log(r)$\\
\hline
\end{tabular}
\label{table_rbf}
\end{center}
\end{table} 
\end{enumerate}

The rest of the article is organized as follows. Section~\ref{section_background} sets the background for the rest of the article, by introducing linear solvers, low-rank matrices and the different hierarchical structures. The next section, Section~\ref{section_key_ideas} presents the key ideas, discusses some previous work on fast direct solvers and highlights the contributions of the present work. A comprehensive discussion on how our approach differs from existing ones is highlighted. Section~\ref{section_algorithm} discusses the extended sparsification \& hierarchical compression algorithm, which is followed by the section on numerical benchmark.

\section{Background}
\label{section_background}
\subsection{Iterative versus direct solvers}
\label{subsection_solvers}
Algorithms for solving linear systems can be broadly classified into: (i) Iterative solvers (ii) Direct solvers. Most of the iterative solvers are based on Krylov subspace techniques~\cite{arnoldi1951principle, hestenes1952methods, saad1986gmres, paige1975solution, van1992bi, freund1991qmr, freund1993transpose} and rely on matrix-vector products. The number of iterations required to achieve a target accuracy is highly problem dependent. In many instances, the large condition number of the matrix or distribution of the eigenvalues of the matrix in the complex plane (e.g., widely spread eigenvalues) result in a large number of iterations. Consequently, pre-conditioners need to be devised to cluster the eigenvalues and thereby accelerate convergence of the iterative solver. Once such a pre-conditioner is in place, the cost of performing the matrix-vector products can be reduced using fast summation techniques like the fast multipole method (FMM)~\cite{beatson1992fast,coifman1993fast,nishimura2002fast}, the Barnes-Hut algorithm~\cite{barnes1986hierarchical}, panel clustering~\cite{hackbusch1989fast}, FFT, wavelet based methods, and others. Of these different fast summation techniques, FMM has often been used in the context of linear systems arising out of boundary integral equations. This is because the Green's function resulting from such integral equations are amenable to FMM~\cite{greengard1987fast, greengard1988rapid, beatson1997short, greengard1997new, cheng1999fast, darve2000fast, darve2001fast, fong2009black, ying2004kernel}. With these in place, fast iterative solvers solve the linear system in linear or almost linear complexity.

Direct solvers on the other hand rely on efficient factorization/elimination \textemdash Gaussian elimination/ LU, QR, etc. \textemdash of the underlying linear system and then using the factorization to obtain the solution. Direct solvers possess a wide range of advantages: (i) The solution can be computed exactly (up-to machine precision), (ii) Scales well with multiple right hand sides (iii) Robust and hence preferred for black-box implementations. However, the major draw back of direct solvers is that it is highly expensive, since the factorization/elimination step scales as $\mathcal{O}(N^3)$ for most of the dense matrices. The focus of this article is to reduce the computational complexity of direct solvers to $\mathcal{O}(N)$ for the class of FMM matrices.

\subsection{FMM matrices / $\mathcal{H}^2$ matrices with strong admissibility (weak low-rank structure)}

\begin{figure}[!htbp]
\centering
\includegraphics[scale=0.1325]{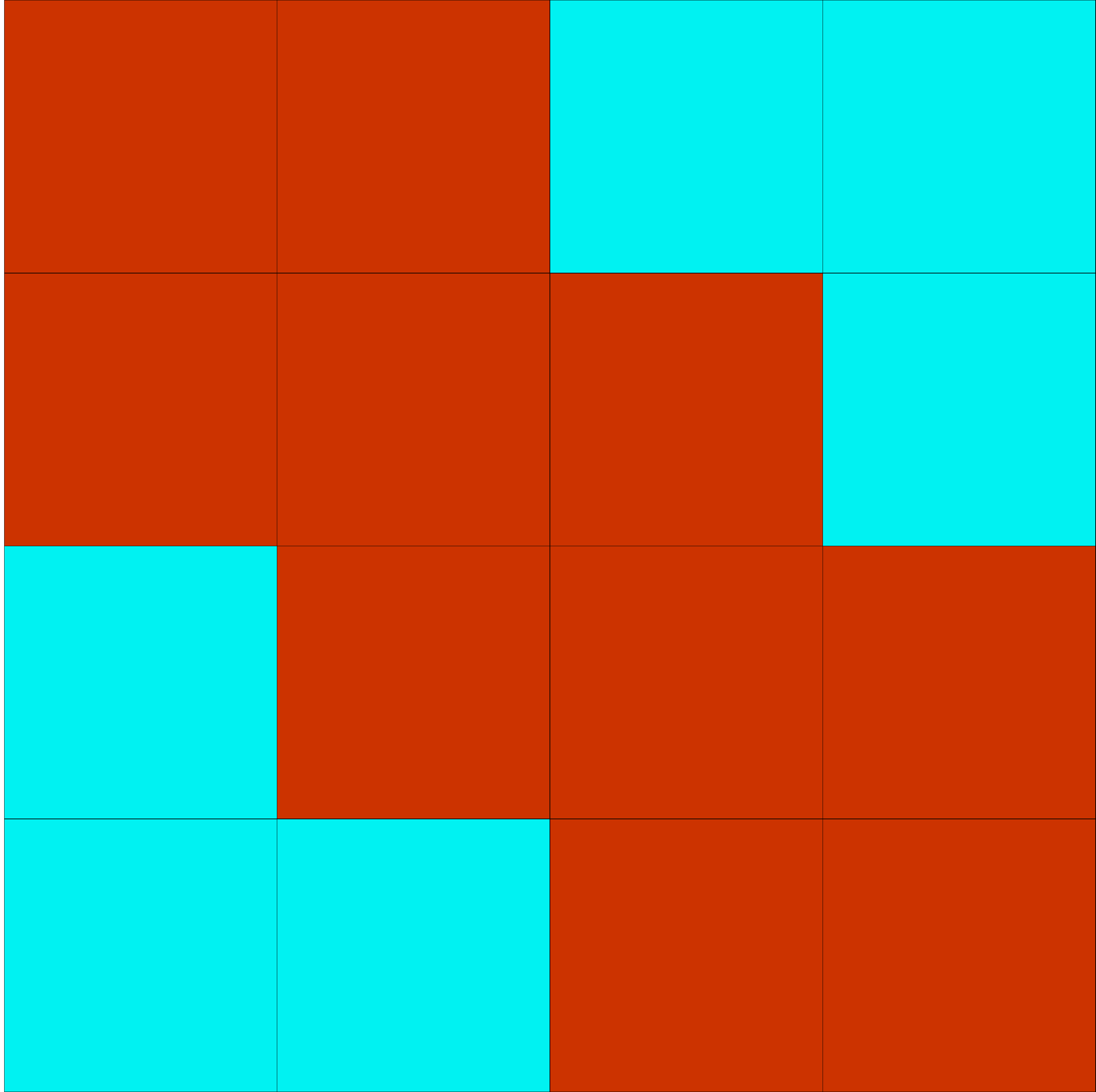}
\includegraphics[scale=0.1325]{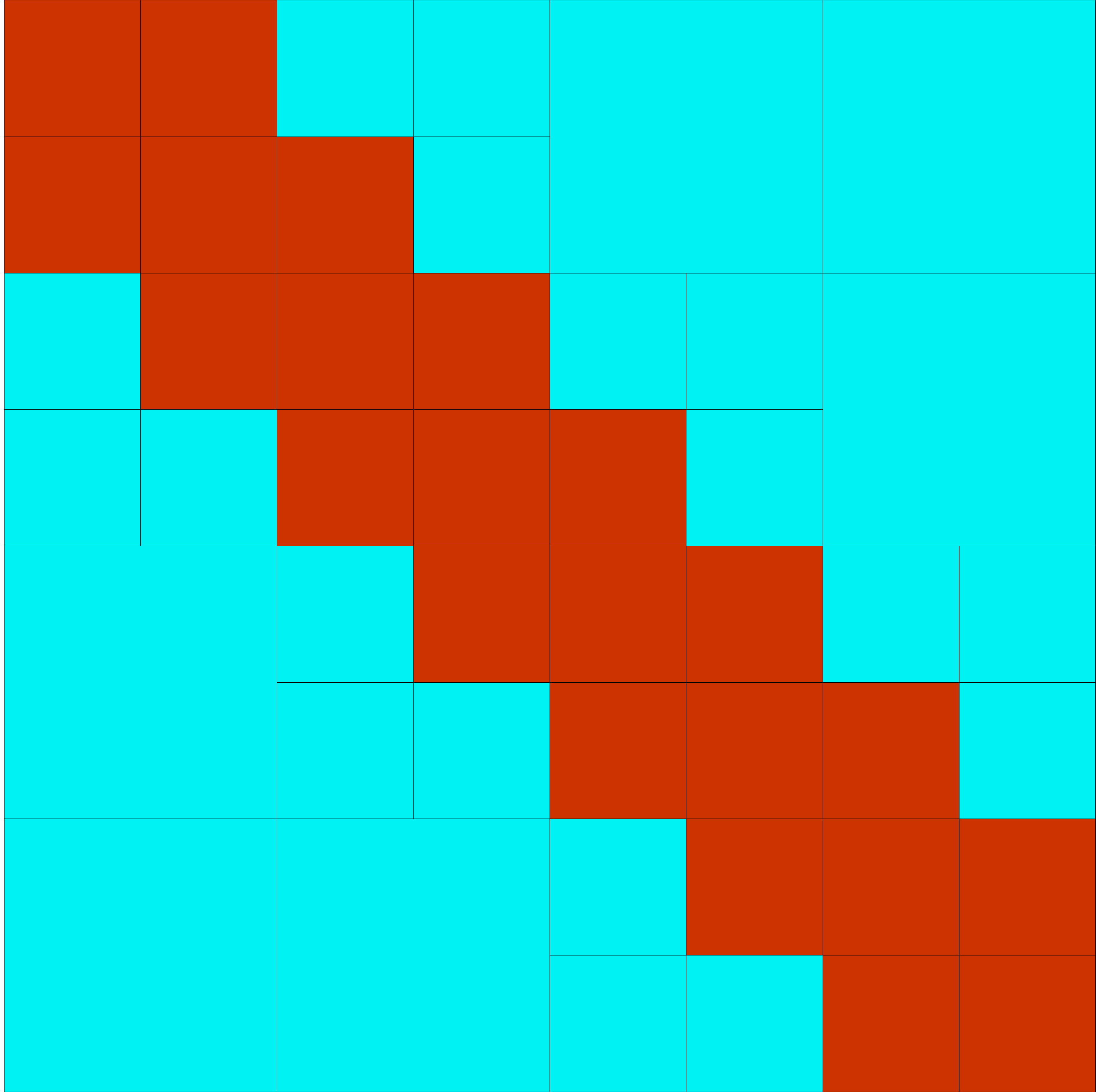}
\includegraphics[scale=0.1325]{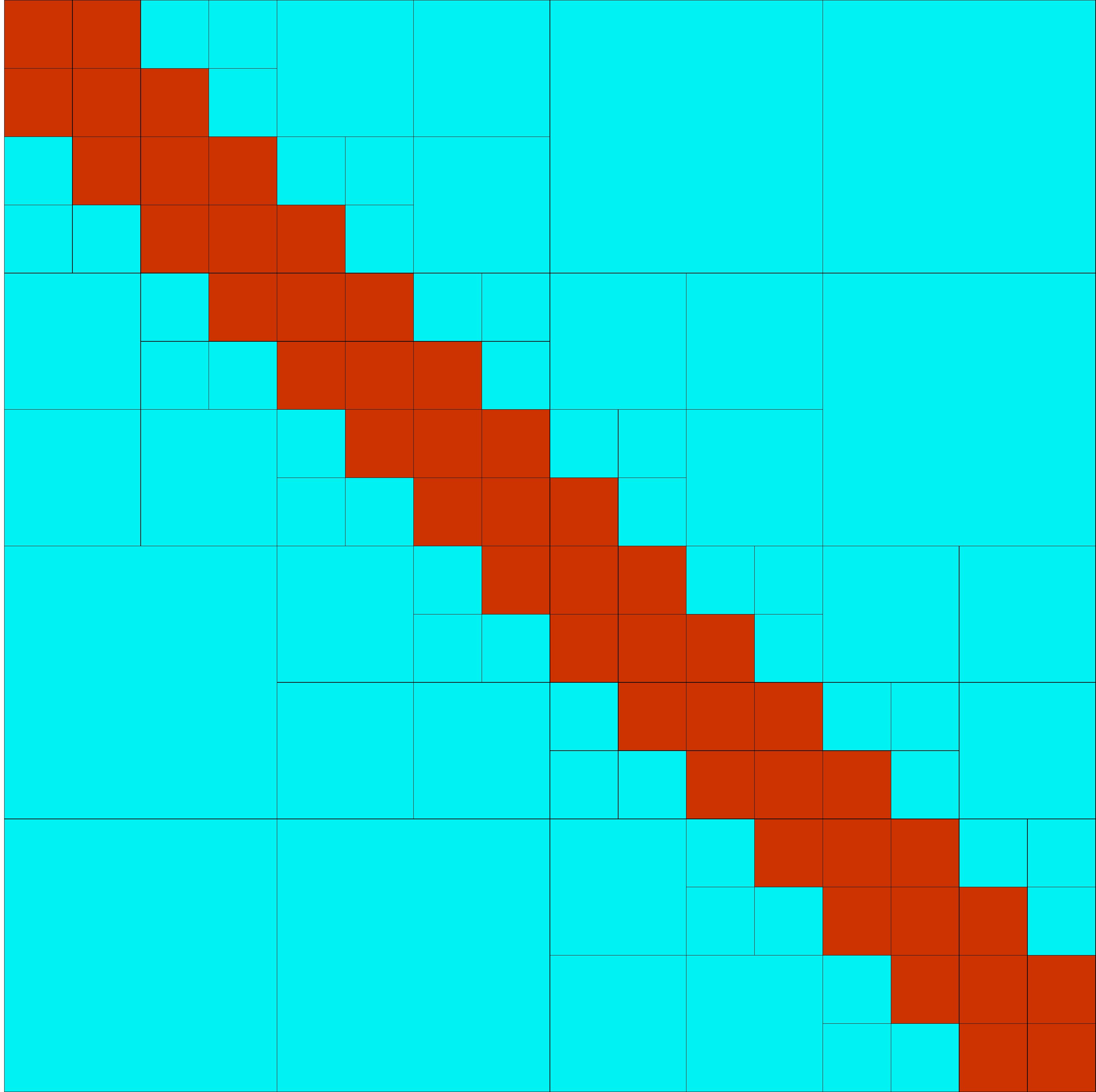}
\tikz{
\draw [fill=fullrank] (0,0) rectangle (0.5,0.5);
\node at (1.5,0.25) {-- Full-rank;};
\draw [fill=lowrank] (3,0) rectangle (3.5,0.5);
\node at (4.5,0.25) {-- Low-rank};
}
\caption{The figure shows the FMM matrix based on a binary tree, where the points lie on an interval with the natural ordering in $1$D. Left: Level $2$; Middle: Level $3$; Right: Level $4$. The figures do not reveal the nested low-rank structure. The row and column basis of the low-rank matrices at a level in the tree can be constructed from the row and column basis of the low-rank matrices of its children.}
\label{figure_1D_H2}
\end{figure}

\begin{figure}[!htbp]
\centering
\includegraphics[scale=0.1]{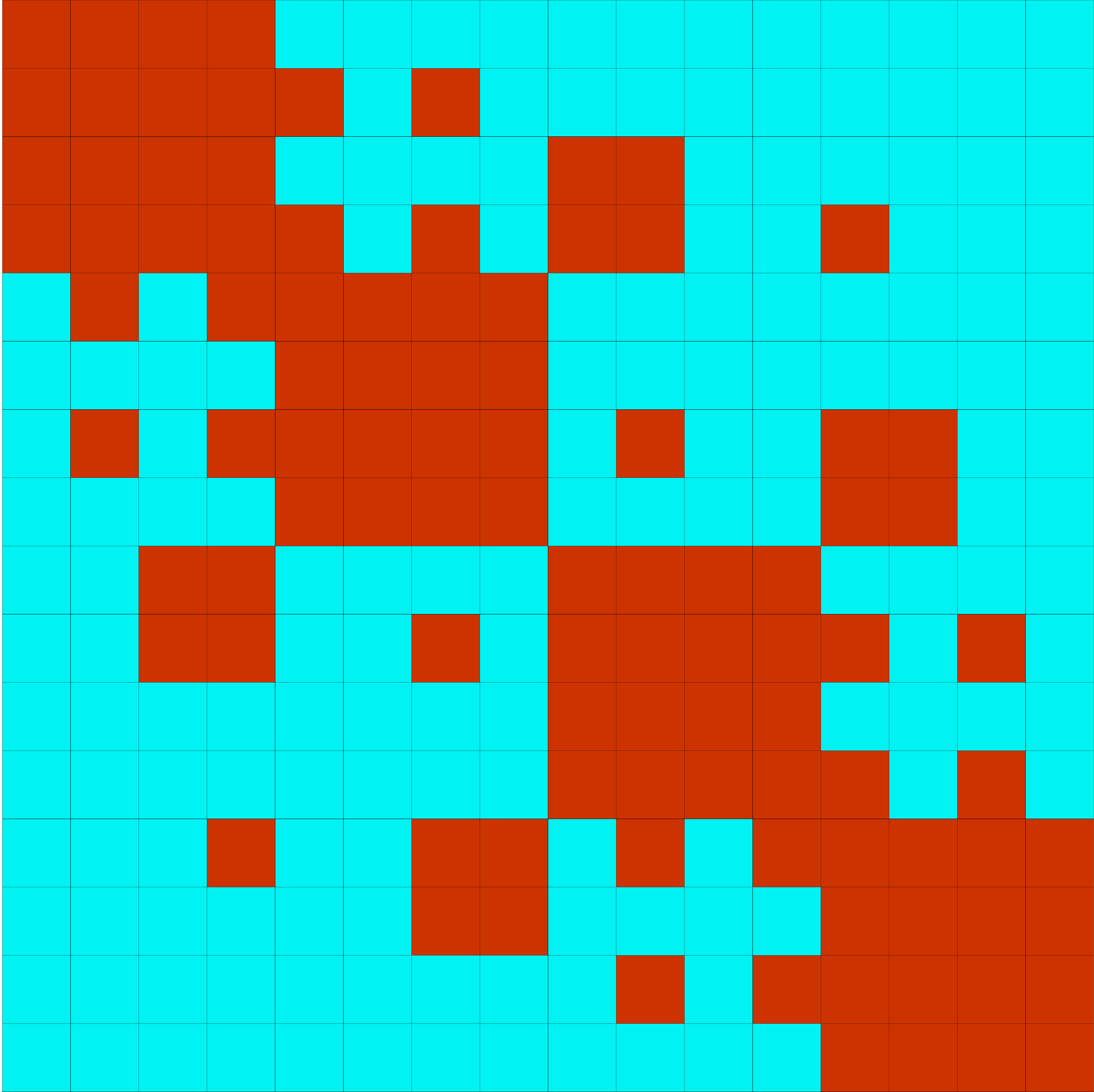}
\includegraphics[scale=0.1]{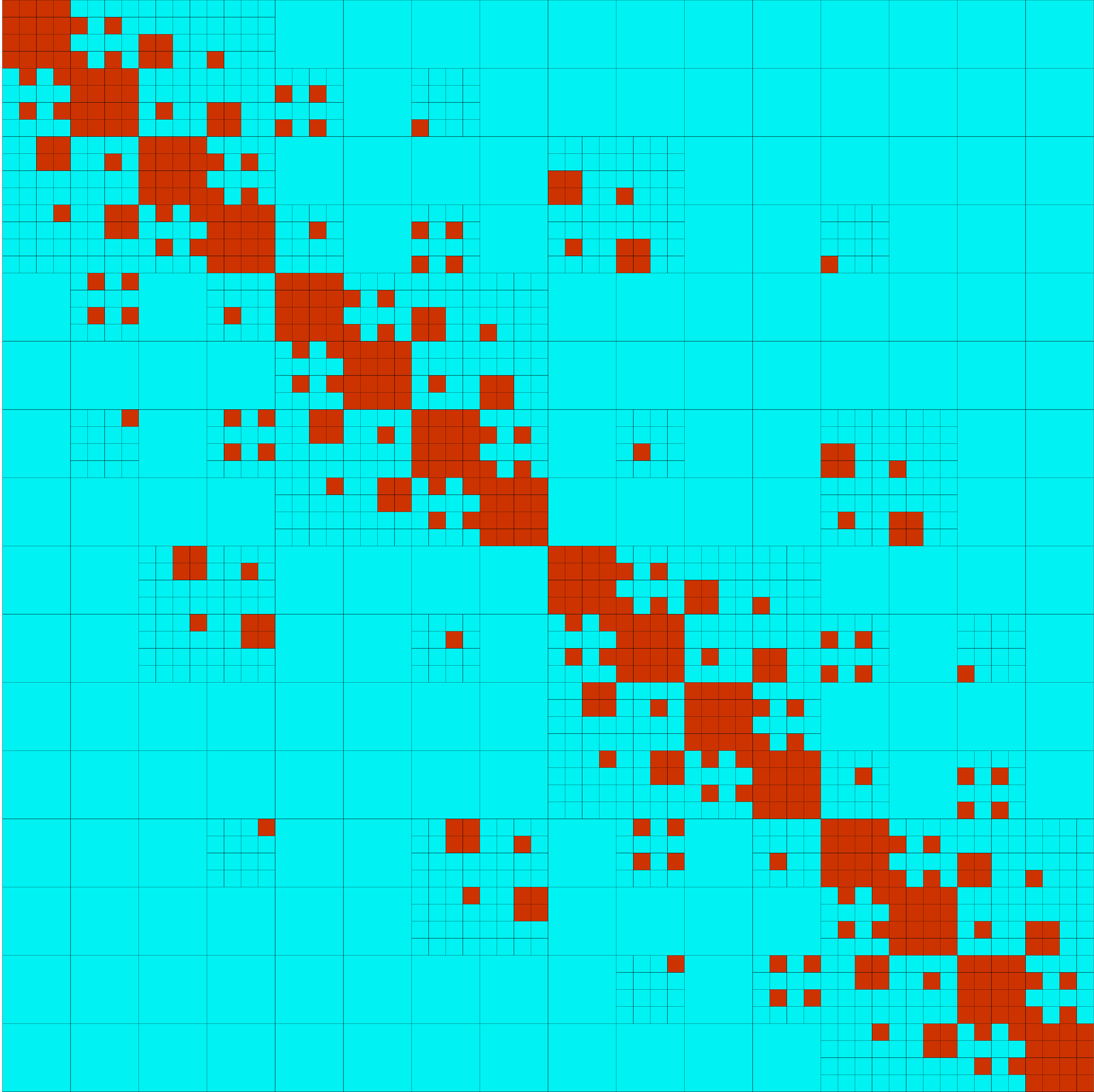}
\tikz{
\draw [fill=fullrank] (0,0) rectangle (0.5,0.5);
\node at (1.5,0.25) {-- Full-rank;};
\draw [fill=lowrank] (3,0) rectangle (3.5,0.5);
\node at (4.5,0.25) {-- Low-rank};
}
\caption{The figure shows the FMM matrix based on a quad tree, where the points lie on a $2$D manifold homeomorphic to a square with Morton ordering/ Z-ordering. Left: Level $2$; Right: Level $3$. The figure on the right does not reveal the nested low-rank structure present in the matrix. The row and column basis of the low-rank matrices at level $2$ in the tree can be constructed from the row and column basis of the low-rank matrices of its children at level $3$.}
\label{figure_2D_H2}
\end{figure}

For readers familiar with FMM~\cite{greengard1987fast}, $\mathcal{H}^2$ matrices are an algebraic version of the matrices encountered in the FMM. A minor difference is that in the FMM, we sub-divide the geometric domain is hierarchically partitioned using a $2^d$ tree in $d$ dimensions. But in the case of $\mathcal{H}^2$ matrices, the tree is not restricted to a $2^d$ tree in $d$ dimensions. Table~\ref{table_1D_H2_2level}, Figure~\ref{figure_1D_H2} and Figure~\ref{figure_2D_H2} illustrate the FMM matrices, on $1$D and $2$D manifolds at different levels in the tree.
\begin{remark}
Note that Figure~\ref{figure_1D_H2} and Figure~\ref{figure_2D_H2} do not explicitly reveal the nested low-rank structure. The low-rank basis at the parent level is constructed from the low-rank basis of its children.
\end{remark}

\begin{figure}[!htbp]
\begin{floatrow}
\ffigbox{%
\includegraphics[scale=0.4]{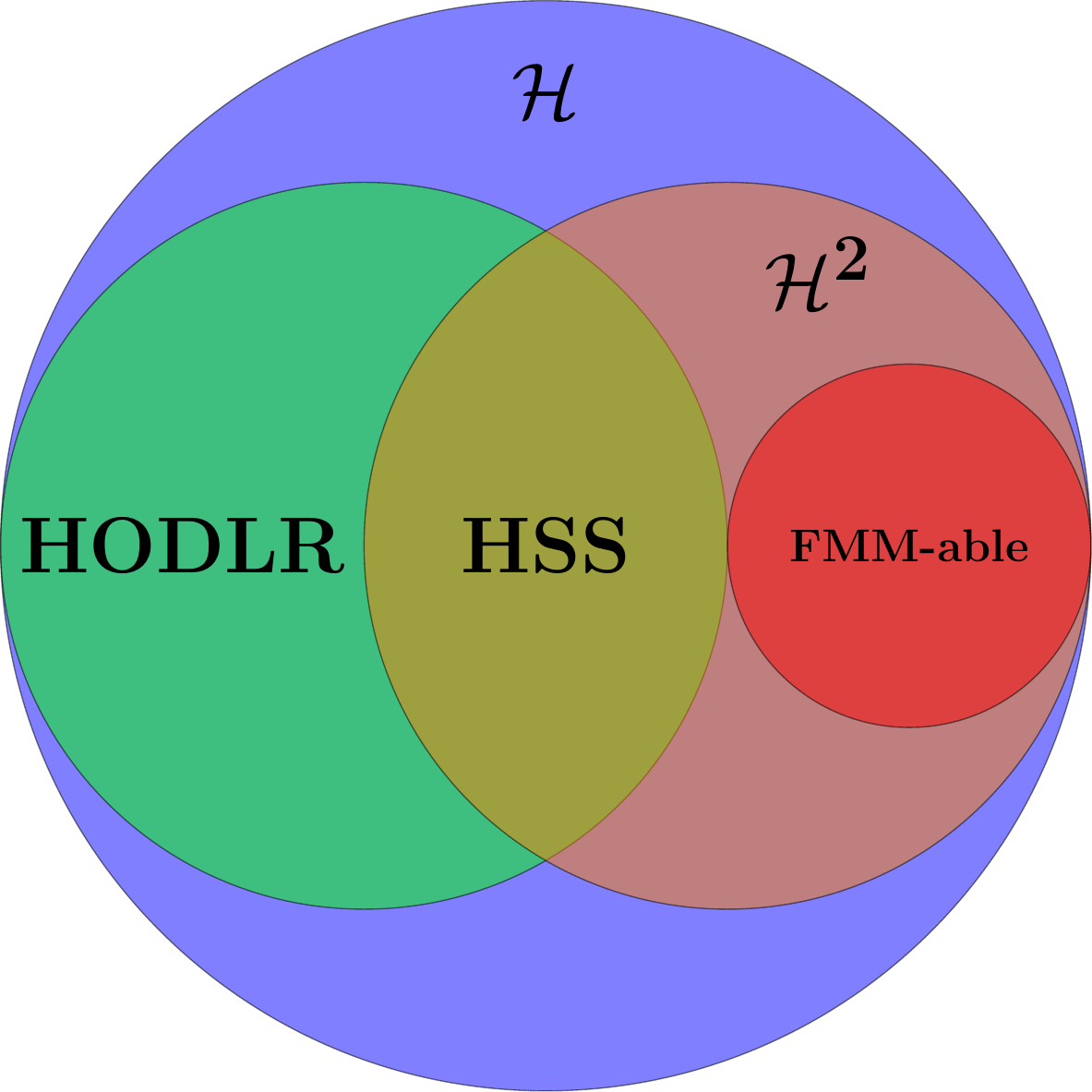}
}{%
\caption{An Euler diagram of the different hierarchical matrices}
\label{figure_Euler_diagram}
}
\capbtabbox{%
\includegraphics[scale=1]{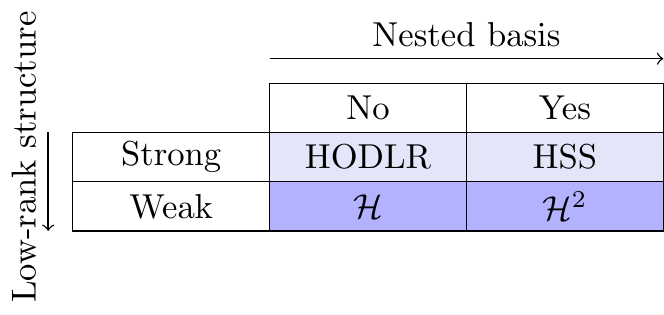}
}{%
\caption{Different hierarchical structures.}
\label{table_hierarchical_structures}
}
\end{floatrow}
\end{figure}

In general, based on the admissibility and nested low-rank structure, hierarchical matrices are classified as shown in Table~\ref{table_hierarchical_structures} and Figure~\ref{figure_Euler_diagram}. For a detailed description of these different hierarchical structures, we refer the readers to Chapter $3$ of~\cite{ambikasaran2013thesis}.
The fast direct solver algorithm presented in the article also extends to other hierarchical structures, though that will not be the focus of the current article. Refer Chapter 7 of~\cite{ambikasaran2013thesis}, where the algorithm is discussed for HSS and HODLR matrices.

\begin{table}[!htbp]
\caption{A $2$ level FMM matrix based on a binary tree, where the points lie on an interval with the natural ordering in $1$D. The matrices $\tilde{K}_{ij}^{(2)}, U_i^{(2)}, V_j^{(2)}$ are low-rank matrices. }
$$
K^{(2)} = \begin{bmatrix}
K_{11}^{(2)} & K_{12}^{(2)} & U_{1}^{(2)} \tilde{K}_{13}^{(2)} V_3^{(2)^T} & U_1^{(2)} \tilde{K}_{14}^{(2)} V_4^{(2)^T}\\
K_{21}^{(2)} & K_{22}^{(2)} & K_{23}^{(2)} & U_2^{(2)} \tilde{K}_{24}^{(2)} V_4^{(2)^T}\\
U_{3}^{(2)} \tilde{K}_{31}^{(2)} V_1^{(2)^T} & K_{32}^{(2)} & K_{33}^{(2)} & K_{34}^{(2)}\\
U_{4}^{(2)} \tilde{K}_{41}^{(2)} V_1^{(2)^T} & U_{4}^{(2)} \tilde{K}_{42}^{(2)} V_2^{(2)^T} & K_{43}^{(2)} & K_{44}^{(2)}\\
\end{bmatrix}
$$
\label{table_1D_H2_2level}
\end{table} 

\section{Key ideas, previous work and contribution of this article}
\label{section_key_ideas}
There are couple of novel ideas exploited in this article to attain linear complexity.
\begin{itemize}
\item
The first key idea is that one form of sparsity can be converted into another form of sparsity, which can then be exploited. To be specific, a $N \times N$ data-sparse FMM matrix can exactly be represented as a $M \times M$ extended structured sparse matrix, where $M \in \mathcal{O}(N)$. The unknowns of the original dense linear system are a subset of the unknowns of the extended sparse linear system. (We use the term sparse matrix in the conventional sense; i.e., a matrix primarily populated with zeros. This should not be confused with data-sparsity, i.e., a matrix with low-rank sub-blocks.)


\begin{figure}[!htbp]
\tikz{
\draw (0,0) -- (9,0);
\draw (0,-0.25) -- (0,0.25);
\draw (3,-0.25) -- (3,0.25);
\draw (6,-0.25) -- (6,0.25);
\draw (9,-0.25) -- (9,0.25);
\node at (1.75,0.25) {Cluster $1$};
\node at (4.75,0.25) {Cluster $2$};
\node at (7.75,0.25) {Cluster $3$};
\draw[-triangle 90] (4.5,-0.125) to[out=-150,in=-30] (2.25,-0.125);
\node at (3.25,-0.75) {$G_{12}$};
\draw[-triangle 90] (7.5,-0.125) to[out=-120,in=-60] (1.5,-0.125);
\node at (5,-1.25) {$G_{13}$};
\draw[-triangle 90] (7.25,-0.125) to[out=210,in=-30] (4.5,-0.125);
\node at (5.75,-0.75) {$G_{23}$};
\draw[-triangle 90] (3.25, 0.125) to[out=60,in=120] (5.75, 0.125);
\node at (4.5,1) {$G_{22}$};
}
\caption{The relevant interactions along a $1$D manifold for Table~\ref{table_numerical_ranks_1D}. $N$ denotes the number of particles in each cluster.}
\end{figure}
\begin{table}[!htbp]
\caption{Ranks of different interactions in $1$D to a precision of $10^{-15}$, i.e., $\dfrac{\sigma_{r+1}}{\sigma_1} < 10^{-15}$, where $r$ is the rank and $\sigma_i$ are the singular values. $a$ was taken as $0.005$ of the length of the interval.}
\begin{tabular}{|c|c|c|c|c|}
\hline
\multirow{2}{*}{Kernels} & \multirow{2}{*}{$N$} & Neighbor & \multicolumn{2}{c|}{Well-separated rank}\\
& & rank & Interaction & Schur complement\\
\hline
 & & $G_{12}$ & $G_{13}$ & $G_{13}$ - $G_{12}G_{22}^{-1} G_{23}$\\
 \hline
 \multirow{3}{*}{Logarithm} & 32 & 14 & 9 & 8\\
 & 64 & 16 & 9 & 9\\
 & 128 & 18 & 9 & 9\\
 \multirow{4}{*}{$\left(\dfrac{r(\log(r)-1)}{a(\log(a)-1)} \right) \chi_{r < a} + \left(\dfrac{\log(r)}{\log(a)} \right)\chi_{t \geq a}$} & 256 & 19 & 8 & 8\\
 & 512 & 21 & 8 & 8\\
 & 1024 & 24 & 8 & 8\\
 & 2048 & 28 & 8 & 8\\
 \hline
 \multirow{3}{*}{Bessel function of second kind} & 32 & 14 & 10 & 9\\
 & 64 & 16 & 10 & 9\\
 & 128 & 18 & 9 & 9\\
 \multirow{4}{*}{$\left(\dfrac{1+y_0(a)}{1+y_0(r)} \right) \chi_{r < a} + \left(\dfrac{y_0(r)}{y_0(a)} \right)\chi_{r \geq a}$} & 256 & 20 & 9 & 9\\
 & 512 & 22 & 9 & 9\\
 & 1024 & 24 & 9 & 9\\
 & 2048 & 29 & 9 & 9\\
 \hline
  \multirow{3}{*}{Inverse quadric} & 32 & 14 & 9 & 9\\
 & 64 & 15 & 9 & 8\\
 & 128 & 14 & 9 & 9\\
 \multirow{4}{*}{$\left(\dfrac{\arctan(r)}{\arctan(a)} \right) \chi_{r < a} + \left(\dfrac{1+a^2}{1+r^2} \right)\chi_{r \geq a}$} & 256 & 14 & 9 & 9\\
 & 512 & 15 & 9 & 8\\
 & 1024 & 18 & 9 & 8\\
 & 2048 & 23 & 9 & 7\\
 \hline
  \multirow{3}{*}{Inverse multi-quadric} & 32 & 14 & 9 & 8\\
 & 64 & 15 & 9 & 9\\
 & 128 & 15 & 9 & 9\\
 \multirow{4}{*}{$\left(\dfrac{\text{arcsinh}(r)}{\text{arcsinh}(a)} \right) \chi_{r < a} + \left(\sqrt{\dfrac{1+a^2}{1+r^2}} \right)\chi_{r \geq a}$} & 256 & 14 & 8 & 8\\
 & 512 & 14 & 8 & 8\\
 & 1024 & 17 & 8 & 8\\
 & 2048 & 22 & 8 & 8\\
 \hline
\end{tabular}
\label{table_numerical_ranks_1D}
\end{table}

\begin{figure}[!htbp]
\tikz[scale=0.6]{
\draw (0,0) rectangle (9,3);
\draw (3,0) -- (3,3);
\draw (6,0) -- (6,3);
\node at (1.75,1.5) {Cluster $1$};
\node at (4.75,1.5) {Cluster $2$};
\node at (7.75,1.5) {Cluster $3$};
\draw[-triangle 90] (4.5,-0.125) to[out=-150,in=-30] (2.25,-0.125);
\node at (3.25,-0.75) {$G_{12}$};
\draw[-triangle 90] (7.5,-0.125) to[out=-120,in=-60] (1.5,-0.125);
\node at (5,-1.25) {$G_{13}$};
\draw[-triangle 90] (7.25,-0.125) to[out=210,in=-30] (4.5,-0.125);
\node at (5.75,-0.75) {$G_{23}$};
\draw[-triangle 90] (3.25, 3.125) to[out=60,in=120] (5.75, 3.125);
\node at (4.5,4) {$G_{22}$};
}
\caption{The relevant interactions along a $2$D manifold for Table~\ref{table_numerical_ranks_2D}. $N$ denotes the number of particles in each cluster.}
\end{figure}

\begin{table}[!htbp]
\caption{Ranks of different interactions in $2$D to a precision of $10^{-15}$, i.e., $\dfrac{\sigma_{r+1}}{\sigma_1} < 10^{-15}$, where $r$ is the rank and $\sigma_i$ are the singular values. $a$ was taken as $0.005$ of the length of the side of the square.}
\begin{tabular}{|c|c|c|c|c|}
\hline
\multirow{2}{*}{Kernels} & \multirow{2}{*}{$N$} & Neighbor & \multicolumn{2}{c|}{Well-separated rank}\\
& & rank & Interaction & Schur complement\\
\hline
 & & $G_{12}$ & $G_{13}$ & $G_{13}$ - $G_{12}G_{22}^{-1} G_{23}$\\
 \hline
 \multirow{3}{*}{Logarithm} & 64 & 26 & 18 & 17\\
 & 121 & 31 & 19 & 19\\
 \multirow{4}{*}{$\left(\dfrac{r(\log(r)-1)}{a(\log(a)-1)} \right) \chi_{r < a} + \left(\dfrac{\log(r)}{\log(a)} \right)\chi_{t \geq a}$} & 225 & 33 & 19 & 19\\
 & 441 & 37 & 19 & 21\\
 & 841 & 41 & 19 & 22\\
 & 1681 & 44 & 18 & 23\\
 & 3249 & 46 & 18 & 24\\
 \hline
 \multirow{3}{*}{Bessel function of second kind} & 64 & 28 & 21 & 19\\
 & 121 & 33 & 21 & 23\\
 & 225 & 37 & 22 & 23\\
 \multirow{4}{*}{$\left(\dfrac{1+y_0(a)}{1+y_0(R)} \right) \chi_{r < a} + \left(\dfrac{y_0(R)}{y_0(a)} \right)\chi_{t \geq a}$} & 441 & 40 & 22 & 25\\
 & 841 & 44 & 22 & 26\\
 & 1681 & 47 & 21 & 27\\
 & 3249 & 49 & 22 & 28\\
 \hline
  \multirow{3}{*}{$3$D Laplace} & 64 & 53 & 36 & 38\\
 & 121 & 68 & 38 & 42\\
 & 225 & 78 & 39 & 45\\
 \multirow{4}{*}{$\left(\dfrac{r}{a} \right) \chi_{r < a} + \left(\dfrac{a}r \right)\chi_{r \geq a}$} & 441 & 88 & 38 & 43\\
 & 841 & 94 & 39 & 45\\
 & 1681 & 103 & 37 & 45\\
 & 3249 & 110 & 37 & 45\\
 \hline
  \multirow{3}{*}{$3$D Helmholtz} & 64 & 56 & 43 & 45\\
 & 121 & 73 & 44 & 42\\
 & 225 & 82 & 45 & 50\\
 \multirow{4}{*}{$\left(\dfrac{r}a \right) \chi_{r < a} + \left({\dfrac{a\cos(r)}{r \cos(a)}} \right)\chi_{r \geq a}$} & 441 & 93 & 45 & 52\\
 & 841 & 101 & 45 & 50\\
 & 1681 & 108 & 45 & 47\\
 & 3249 & 116 & 43 & 51\\
 \hline
  \multirow{3}{*}{Biharmonic} & 64 & 38 & 25 & 21\\
 & 121 & 43 & 25 & 20\\
 & 225 & 49 & 27 & 20\\
 \multirow{4}{*}{$\left(\dfrac{r^3(3\log(r)-1)}{a^3(3\log(a)-1)} \right) \chi_{r < a} + \left({\dfrac{r^2\log(r)}{a^2\log(a)}} \right)\chi_{r \geq a}$} & 441 & 52 & 26 & 20\\
 & 841 & 57 & 27 & 20\\
 & 1681 & 58 & 25 & 20\\
 & 3249 & 62 & 25 & 20\\
 \hline
\end{tabular}
\label{table_numerical_ranks_2D}
\end{table} 

\item
The second \emph{key ingredient} is that, when performing the elimination of unknowns in the extended sparse linear system, the interaction between the unknowns corresponding to the well-separated clusters at all stages in the elimination process can be efficiently compressed as low-rank, which is validated in Tables~\ref{table_numerical_ranks_1D} and~\ref{table_numerical_ranks_2D} for different kernel functions.

This implies, after an appropriate ordering of equations and unknowns, while we perform the elimination, the fill-in that occurs in the elimination process corresponding to well-separated clusters can be compressed and efficiently represented as a low-rank matrix. As shown in Section~\ref{section_algorithm}, the ordering of the equations and the unknowns in the extended sparser matrix is strongly related to the local and multipole coefficients in the fast multipole method, which is different from the one obtained using the nested dissection approach.
\end{itemize}

Before discussing our algorithm, we present a brief discussion of previous works in this direction and our new contribution. The idea of extended sparsification has been considered before in the article by Chandrasekaran et al.~\cite{chandrasekaran2006fast}, though only in the context of HSS matrices, which are strict sub-class of $\mathcal{H}^2$ matrices, i.e., has the additional constraint that the interaction between all (not just the ``well-separated clusters") are low-rank. The algorithm for HSS matrices is fairly easier since in the elimination process there are \emph{no new fill-ins}. However, hierarchically semi-separable matrices are restricted to one-dimensional applications. In our approach, we deal with the larger class of FMM matrices, which model a large class of hierarchical matrices in all dimensions.

Further it is non trivial to extend the algorithm presented in Chandrasekaran et al.~\cite{chandrasekaran2006fast} to FMM matrices. To achieve linear complexity for FMM matrices, \emph{it is to be emphasized that the second step mentioned in Section~\ref{section_key_ideas} is highly crucial.} For instance, Pals, in his thesis~\cite{pals2004multipole}, follows a similar approach of representing the matrix arising out of fast multipole method as an extended sparser matrix, but does not exploit the fact that the fill-ins can be compressed as we proceed through the algorithm. Instead Pals~\cite{pals2004multipole} rely on nested dissection using METIS to solve the sparse linear system and show that even though the extended sparser approach is faster, the scaling of the algorithm is still in fact very close to $\mathcal{O}(N^3)$ (Refer Chapter $4$ of~\cite{pals2004multipole}). This is due to the fact that the conventional nested dissection approach doesn't exploit the fact that the fill-ins are low-rank.

Greengard et al.~\cite{greengard2009fast} present the idea of representing the dense matrix as an extended sparse matrix. The article presents a single-level fast solver whose scaling is not $\mathcal{O}(N)$ and due to which compressing the fill-ins, which is important to extend the strategy in a multi-level setting, is not discussed.

Ho and Greengard~\cite{ho2012fast} discuss the use of extended sparsification technique mentioned in~\cite{chandrasekaran2006fast1}, i.e., based on HSS representations and a recursive skeletonization approach, in the context of integral equations, where the compression is obtained using the interpolative decomposition technique. The computational complexity scales like $\mathcal{O}(N)$ on $1$D manifolds.

Ambikasaran's thesis~\cite{ambikasaran2013thesis} discusses the $\mathcal{O}(N)$ extended sparsification technique for HSS matrices and extends it to the bigger class of HODLR matrices at a computational complexity of $\mathcal{O}(N \log^2 N)$, with application to interpolation using radial basis functions.

There have also been other attempts not based on the extended sparsification technique. Broadly speaking, these techniques rely on factorizing the matrix instead of introducing additional variables and we refer the readers to the work by Hackbusch and coworkers~\cite{hackbusch2000sparse,hackbusch2000sparse1,hackbusch2002data}. Ambikasaran \& Darve~\cite{ambikasaran2013fast} and Kong et al.~\cite{kong2011adaptive} discuss an $\mathcal{O}(N \log^2 N)$ algorithm on $1$D manifolds in the context of radial basis function interpolation and integral equations respectively using HODLR matrices. The algorithm is based on the Sherman-Morrison-Woodbury formula~\cite{woodbury1950inverting,hager1989updating}. Ambikasaran \& Darve~\cite{ambikasaran2013fast} also extend the algorithm to the context of HSS matrices at a computational complexity of $\mathcal{O}(N \log N)$.

A $ULV^T$ decomposition ($U$ and $V$ are unitary matrices, and $L$ is lower triangular) of a HSS matrix is discussed in work by Chandrasekaran et al.~\cite{chandrasekaran2006fast}. The key ingredient of their algorithm is to recognize that with a low-rank approximation of the form $U B V^H$, it is possible to apply a unitary transformation to have the last set of rows of $U$ to be non-zero. This is then applied recursively in a ``bottom-up" fashion to attain the factorization at a computational complexity of $\mathcal{O}(N)$.

The work by Rokhlin and Martinsson~\cite{martinsson2005fast} constructs an $\mathcal{O}(N)$ fast direct solver for boundary integral equations in two-dimensions (i.e., one dimensional manifold) making use of off-diagonal low-rank blocks, i.e., a HSS representation of the matrix. The algorithm constructs the inverse using a compressed block factorization that takes advantage of the low-rank off-diagonal blocks to factor the matrix.

To summarize, most of the previous work on fast direct solvers~\cite{ambikasaran2013fast,kong2011adaptive, ho2012fast,martinsson2005fast,chandrasekaran2006fast,chandrasekaran2006fast1,chandrasekaran2002fast} rely on HSS approach. The main drawback of the HSS based fast direct solver is that it is restrictive, especially for applications involving dense matrices arising in $2$D and $3$D. In particular, the rank of the off-diagonal blocks grow as $\mathcal{O}(N^{1/2})$ and $\mathcal{O}(N^{2/3})$ in $2$D and $3$D respectively.

The main contribution of the algorithm discussed in this paper is that we abandon the HSS matrix framework and work with the more general class of FMM matrices throughout the algorithm. There has been some previous work on $\mathcal{H}^2$ matrices~\cite{hackbusch2002h2,borm2010efficient,borm2006matrix}, but they are mainly restricted to almost linear complexity matrix-vector products and fast iterative solvers. To our knowledge, this is the first direct solver for FMM matrices in all dimensions. It is to be noted that B\"{o}rm~\cite{borm2006matrix} discusses an algorithm for $\mathcal{H}^2$ matrix-matrix multiplications. Using this a fast direct solver for $\mathcal{H}^2$ matrices could be constructed but the pre-factor tends to be large, since the $\mathcal{H}^2$ matrix-matrix multiplications are expensive, despite their linear complexity. The numerical benchmarks of our algorithm indicate that the pre-factor in the scaling is not that large and large problems can be solved in a reasonable amount of time.

\begin{remark}
The most important aspect of the algorithm discussed in this article is that \textbf{at all stages} in the algorithm, we only represent the interaction between well-separated clusters as ``low-rank", i.e., we only rely on compressing the interaction between ``well-separated" clusters. Hence, the ``low-rank" matrices considered in our algorithm are ``truly" low-rank, i.e., the rank of these matrices is independent of $N$, the cluster size.
\end{remark}

\section{Fast direct solver for FMM matrices}
\label{section_algorithm}
In this section, we first look at how the dense FMM matrix can be interpreted as an extended sparser matrix. It is worth recalling that in the FMM tree data structure in $d$ dimensions, we have a $2^d$ tree and there are local and multipole coefficients at each level in the tree. To form the extended sparser system, we introduce these local and multipole coefficients as unknowns, and the corresponding set of relations between them as equations. \emph{We then present a new ordering of the equations/relations, which is different from the nested dissection ordering for sparse linear systems.}
\begin{remark}
We first choose to explain the algorithm in $1$D and using pictures of matrices and their corresponding graphs. This is done for a couple of reasons.
\begin{itemize}
	\item
	Explaining in $1$D succinctly captures almost all the key features of the algorithm (Refer remark~\ref{remark_1Ddiff}).
	\item
	Pictures of matrices and their graphs provides an easy way to internalize the algorithm.
\end{itemize}
The general algorithm applicable in any dimension is presented later in section~\ref{subsection_Algorithm}.
\end{remark}

\subsection{Illustration in $1$D}
Consider the linear equation~\eqref{equation_1D_H2_2level} obtained from a $2$ level FMM matrix on an interval in $1$D.

\subsubsection{Ordering of equations and unknowns}
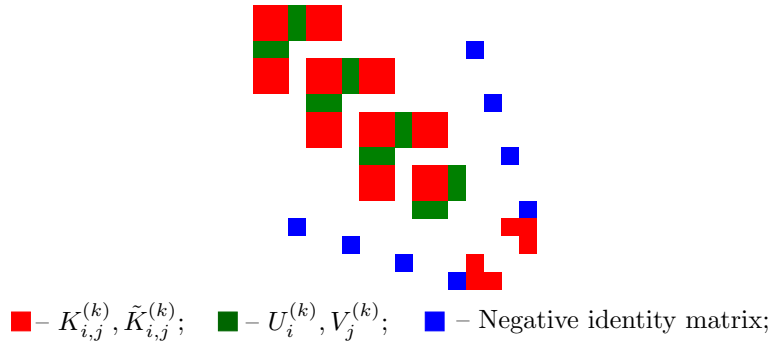
\begin{figure}[!htbp]
\centering
\includegraphics[scale=0.15]{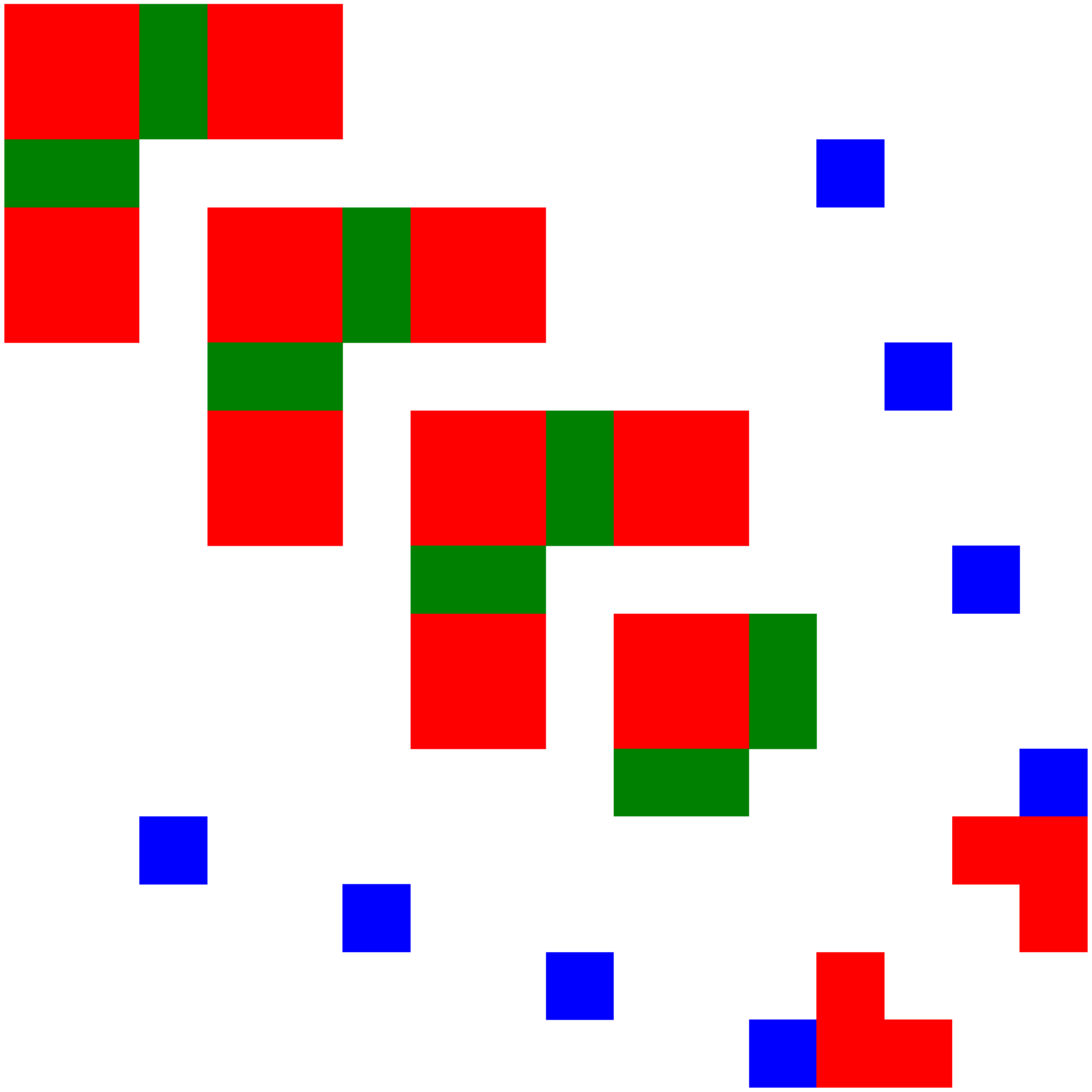}\\
\tikz{
\draw [draw=red,fill=red] (0,0) rectangle (0.25,0.25);
\node at (1.3125,0.125) {-- $K_{i,j}^{(k)}, \tilde{K}_{i,j}^{(k)}$;};
\draw [draw=darkgreen,fill=darkgreen] (2.75,0) rectangle (3,0.25);
\node at (4.0625,0.125) {-- $U_{i}^{(k)}, V_j^{(k)}$;};
\draw [draw=blue,fill=blue] (5.5,0) rectangle (5.75,0.25);
\node at (8,0.125) {-- Negative identity matrix;};
}
\caption{Sparsity pattern of FMM matrix at level $3$ arising out of a $1$D manifold homeomorphic to an interval represented as an extended sparse matrix after appropriate ordering of equations and unknowns.}
\label{figure_1D_H2_2level_extended_sparsity_reordered}
\end{figure}

\begin{align}
\begin{bmatrix}
K_{11}^{(2)} & K_{12}^{(2)} & U_{1}^{(2)} \tilde{K}_{13}^{(2)} V_3^{(2)^T} & U_1^{(2)} \tilde{K}_{14}^{(2)} V_4^{(2)^T}\\
K_{21}^{(2)} & K_{22}^{(2)} & K_{23}^{(2)} & U_2^{(2)} \tilde{K}_{24}^{(2)} V_4^{(2)^T}\\
U_{3}^{(2)} \tilde{K}_{31}^{(2)} V_1^{(2)^T} & K_{32}^{(2)} & K_{33}^{(2)} & K_{34}^{(2)}\\
U_{4}^{(2)} \tilde{K}_{41}^{(2)} V_1^{(2)^T} & U_{4}^{(2)} \tilde{K}_{42}^{(2)} V_2^{(2)^T} & K_{43}^{(2)} & K_{44}^{(2)}\\
\end{bmatrix}
\begin{bmatrix}
x_1^{(2)}\\
x_2^{(2)}\\
x_3^{(2)}\\
x_4^{(2)}
\end{bmatrix}
=
\begin{bmatrix}
b_1^{(2)}\\
b_2^{(2)}\\
b_3^{(2)}\\
b_4^{(2)}
\end{bmatrix}
\label{equation_1D_H2_2level}
\end{align}
Let us now introduce the multipoles and locals for each cluster. The multipoles and locals for each cluster are given in Equation~\eqref{equation_multipoles_1D_H2_2level} and Equation~\eqref{equation_locals_1D_H2_2level}.
\begin{align}
y_1^{(2)} & = V_1^{(2)^T}x_1; y_2^{(2)} = V_2^{(2)^T}x_2; y_3^{(2)} = V_3^{(2)^T}x_3; y_4^{(2)} = V_4^{(2)^T}x_4
\label{equation_multipoles_1D_H2_2level}\\
z_1^{(2)} & = \tilde{K}_{13}^{(2)}y_3^{(2)} + \tilde{K}_{14}^{(2)}y_4^{(2)}; z_2^{(2)} = \tilde{K}_{24}^{(2)}y_4^{(2)}; z_3^{(2)} = \tilde{K}_{31}^{(2)}y_1^{(2)}; z_4^{(2)} = \tilde{K}_{41}^{(2)}y_1^{(2)} + \tilde{K}_{32}^{(2)}y_2^{(2)}
\label{equation_locals_1D_H2_2level}
\end{align}
Introducing Equations~\eqref{equation_multipoles_1D_H2_2level} \&~\eqref{equation_locals_1D_H2_2level} in Equation~\eqref{equation_1D_H2_2level}, gives us Equation~\eqref{equation_1D_H2_2level_extended_sparsity}.
\begin{equation}
\begin{bmatrix}
K_{11}^{(2)} & K_{12}^{(2)} & 0 & 0 & U_1^{(2)} & 0 & 0 & 0 & 0 & 0 & 0 & 0 \\
K_{21}^{(2)} & K_{22}^{(2)} & K_{23}^{(2)} & 0 & 0 & U_2^{(2)} & 0 & 0 & 0 & 0 & 0 & 0 \\
0 & K_{32}^{(2)} & K_{33}^{(2)} & K_{34}^{(2)} & 0 & 0 & U_3^{(2)} & 0 & 0 & 0 & 0 & 0 \\
0 & 0 & K_{43}^{(2)} & K_{44}^{(2)} & 0 & 0 & 0 & U_4^{(2)} & 0 & 0 & 0 & 0 \\
0 & 0 & 0 & 0 & -I & 0 & 0 & 0 & 0 & 0 & \tilde{K}_{13}^{(2)} & \tilde{K}_{14}^{(2)}\\
0 & 0 & 0 & 0 & 0 & -I & 0 & 0 & 0 & 0 & 0 & \tilde{K}_{24}^{(2)}\\
0 & 0 & 0 & 0 & 0 & 0 & -I & 0 & \tilde{K}_{31}^{(2)} & 0 & 0 & 0\\
0 & 0 & 0 & 0 & 0 & 0 & 0 & -I & \tilde{K}_{41}^{(2)} & \tilde{K}_{42}^{(2)} & 0 & 0\\
V_1^{(2)^T} & 0 & 0 & 0 & 0 & 0 & 0 & 0 & -I & 0 & 0 & 0\\
0 & V_2^{(2)^T} & 0 & 0 & 0 & 0 & 0 & 0 & 0 & -I & 0 & 0\\
0 & 0 & V_3^{(2)^T} & 0 & 0 & 0 & 0 & 0 & 0 & 0 & -I & 0\\
0 & 0 & 0 & V_4^{(2)^T} & 0 & 0 & 0 & 0 & 0 & 0 & 0 & -I\\
\end{bmatrix}
\begin{bmatrix}
x_1^{(2)}\\
x_2^{(2)}\\
x_3^{(2)}\\
x_4^{(2)}\\
z_1^{(2)}\\
z_2^{(2)}\\
z_3^{(2)}\\
z_4^{(2)}\\
y_1^{(2)}\\
y_2^{(2)}\\
y_3^{(2)}\\
y_4^{(2)}
\end{bmatrix}
=
\begin{bmatrix}
b_1^{(2)}\\
b_2^{(2)}\\
b_3^{(2)}\\
b_4^{(2)}\\
0\\
0\\
0\\
0\\
0\\
0\\
0\\
0
\end{bmatrix}
\label{equation_1D_H2_2level_extended_sparsity}
\end{equation}
Note that with the default ordering of equations and unknowns the matrix in Equation~\eqref{equation_1D_H2_2level_extended_sparsity} is asymmetric in terms of fill-in. We now reorder the equations and unknowns to make the matrix into a symmetric matrix as shown in Equation~\eqref{equation_1D_H2_2level_extended_sparsity_reordered}.
\begin{equation}
\begin{bmatrix}
K_{11}^{(2)} & U_1^{(2)} & K_{12}^{(2)} & 0 & 0 & 0 & 0 & 0 & 0 & 0 & 0 & 0 \\
V_1^{(2)^T} & 0 & 0 & 0 & 0 & 0 & 0 & 0 & -I & 0 & 0 & 0\\
K_{21}^{(2)} & 0 & K_{22}^{(2)} & U_2^{(2)} & K_{23}^{(2)} & 0 & 0 & 0 & 0 & 0 & 0 & 0\\
0 & 0 & V_2^{(2)^T} & 0 & 0 & 0 & 0 & 0 & 0 & -I & 0 & 0\\
0 & 0 & K_{32}^{(2)} & 0 & K_{33}^{(2)} & U_3^{(2)} & K_{34}^{(2)} & 0 & 0 & 0 & 0 & 0\\
0 & 0 & 0 & 0 & V_3^{(2)^T} & 0 & 0 & 0 & 0 & 0 & -I & 0\\
0 & 0 & 0 & 0 & K_{43}^{(2)} & 0 & K_{44}^{(2)} & U_4^{(2)} & 0 & 0 & 0 & 0 \\
0 & 0 & 0 & 0 & 0 & 0 & V_4^{(2)^T} & 0 & 0 & 0 & 0 & -I\\
0 & -I & 0 & 0 & 0 & 0 & 0 & 0 & 0 & 0 & \tilde{K}_{13}^{(2)} & \tilde{K}_{14}^{(2)}\\
0 & 0 & 0 & -I & 0 & 0 & 0 & 0 & 0 & 0 & 0 & \tilde{K}_{24}^{(2)}\\
0 & 0 & 0 & 0 & 0 & -I & 0 & 0 & \tilde{K}_{31}^{(2)} & 0 & 0 & 0\\
0 & 0 & 0 & 0 & 0 & 0 & 0 & -I & \tilde{K}_{41}^{(2)} & \tilde{K}_{42}^{(2)} & 0 & 0\\
\end{bmatrix}
\begin{bmatrix}
x_1^{(2)}\\
z_1^{(2)}\\
x_2^{(2)}\\
z_2^{(2)}\\
x_3^{(2)}\\
z_3^{(2)}\\
x_4^{(2)}\\
z_4^{(2)}\\
y_1^{(2)}\\
y_2^{(2)}\\
y_3^{(2)}\\
y_4^{(2)}
\end{bmatrix}
=
\begin{bmatrix}
b_1^{(2)}\\
0\\
b_2^{(2)}\\
0\\
b_3^{(2)}\\
0\\
b_4^{(2)}\\
0\\
0\\
0\\
0\\
0
\end{bmatrix}
\label{equation_1D_H2_2level_extended_sparsity_reordered}
\end{equation}
The extended sparse matrix in the above equation is pictorially represented as shown in Figure~\ref{figure_1D_H2_2level_extended_sparsity_reordered}.
The color code as shown in Figure~\ref{figure_1D_H2_2level_extended_sparsity_reordered} will be followed in the rest of the article as well; Red in the extended sparse matrix denotes the matrices $K_{i,j}^{(k)}$ \& $\tilde{K}_{i,j}^{(k)}$, i.e., the direct interactions and M2L operators; Dark green denotes the interpolation/scatter/L2L/L2P operator \& the anterpolation/gather/M2M/P2M operator; Blue denotes the negative identity matrix.
\begin{figure}[!htbp]
\centering
\includegraphics[scale=0.5]{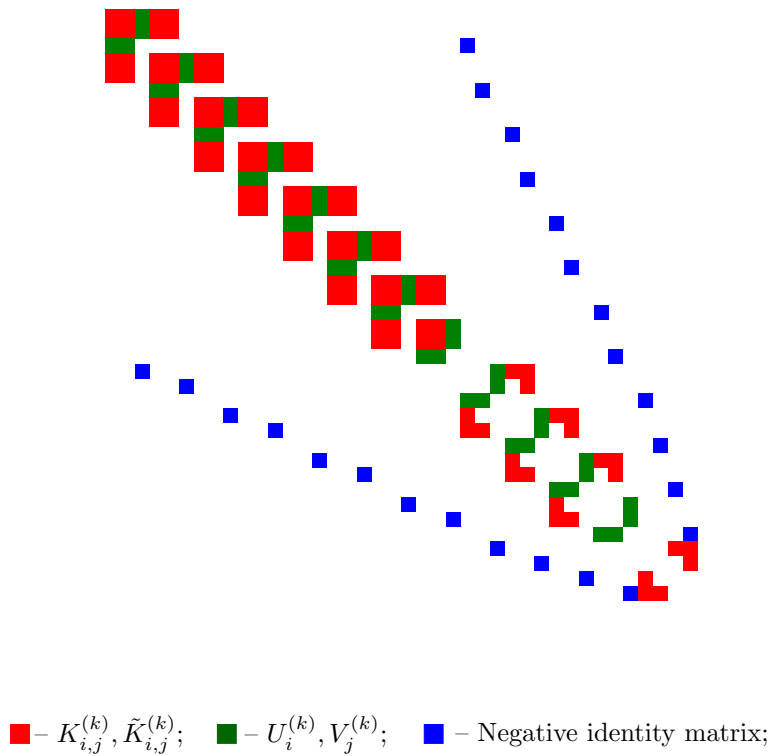}\\
\tikz{
\draw [draw=red,fill=red] (0,0) rectangle (0.25,0.25);
\node at (1.3125,0.125) {-- $K_{i,j}^{(k)}, \tilde{K}_{i,j}^{(k)}$;};
\draw [draw=darkgreen,fill=darkgreen] (2.75,0) rectangle (3,0.25);
\node at (4.0625,0.125) {-- $U_{i}^{(k)}, V_j^{(k)}$;};
\draw [draw=blue,fill=blue] (5.5,0) rectangle (5.75,0.25);
\node at (8,0.125) {-- Negative identity matrix;};
}
\caption{Sparsity pattern of FMM matrix at level $3$ arising out of a $1$D manifold homeomorphic to an interval represented as an extended sparse matrix after appropriate ordering of equations and unknowns.}
\label{figure_1D_H2_3level_extended_sparsity_reordered}
\end{figure} 

\subsubsection{Elimination and hierarchical compression}

Before we describe the algorithm, we would like to mention a key fact.
\begin{figure}[!htbp]
\centering
{\subfigure[Pattern of FMM matrix at level $2$ arising out elimination of leaf level without any compression of the fill-in of the matrix in Figure~\ref{figure_1D_H2_3level_extended_sparsity_reordered}]{
\tikz{
\node at (2.5,3) {\includegraphics[scale=0.2]{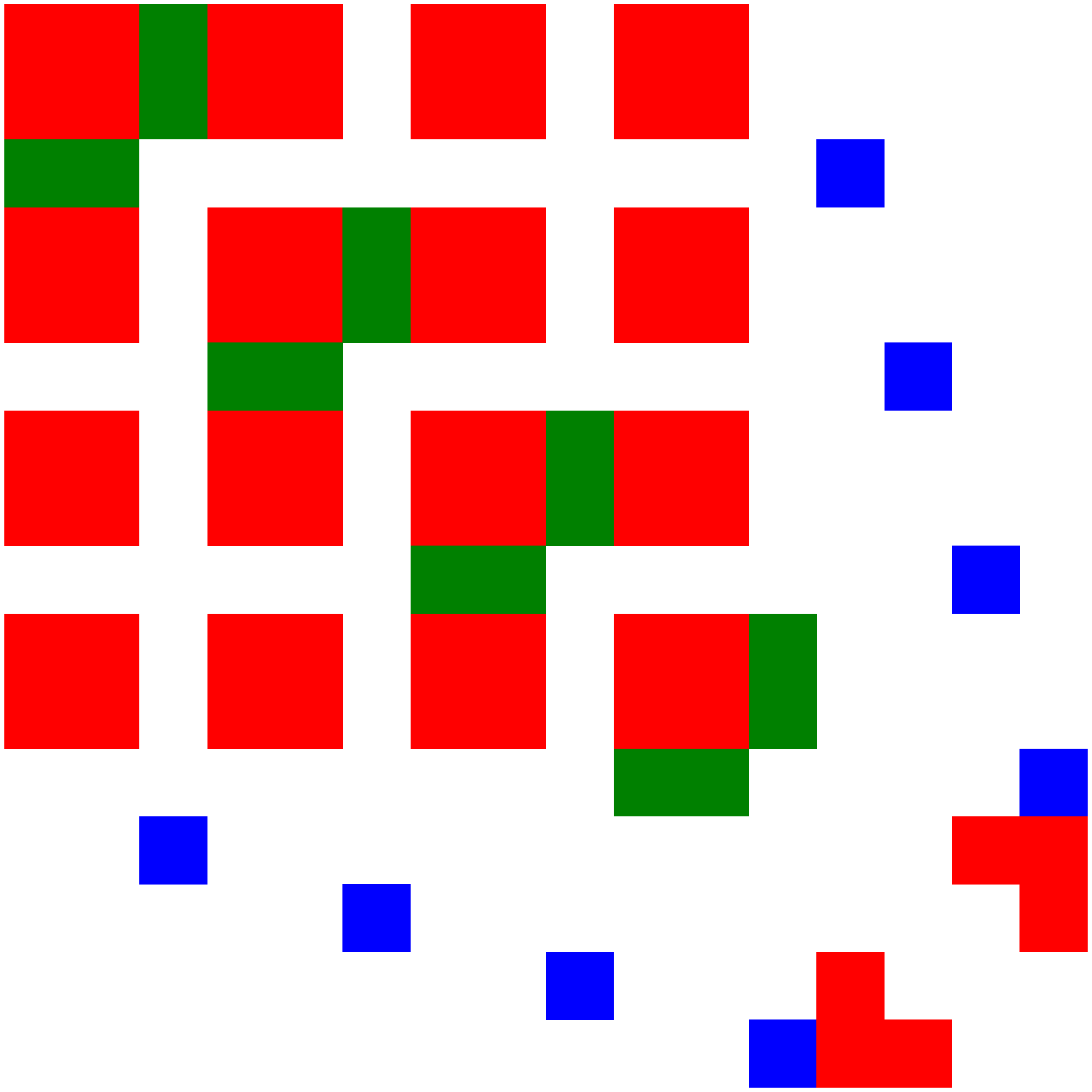}};
}
\label{figure_1D_H2_2level_fillin_extended_sparsity_reordered}
}
}
{\subfigure[Fill-ins that are detrimental to the scaling of the algorithm for the matrix in Figure~\ref{figure_1D_H2_2level_fillin_extended_sparsity_reordered} are highlighted. The key feature of the algorithm is that these fill-ins are avoided by compressing appropriate fill-ins corresponding to well-separated interactions, as we proceed with the elimination.]{
\tikz{
\node at (5,3) {\includegraphics[scale=0.2]{./images/1D_H2_3level_fillin.pdf}};
\fill [draw=blue, fill=blue, fill opacity=0.2, rounded corners] (0.875+3.375,4.75) -- (5.125,4.75) -- (5.125,3.75) -- (6.125,3.75) -- (6.125,5.625) -- (0.875+3.375, 5.625) -- cycle;
\fill [draw=blue, fill=blue, fill opacity=0.2, rounded corners] (2.25,3.75) -- (2.25,1.875) -- (4.25,1.875) -- (4.25,2.875) -- (3.25,2.875) -- (3.25,3.75) -- cycle;
}
}
}
\caption{Extended sparse matrix at level $2$ after eliminating the nodes of the matrix in Figure~\ref{figure_1D_H2_3level_extended_sparsity_reordered} corresponding to level $3$.}
\end{figure} 

%

\begin{remark}
It is important to note that if we feed in the sparse matrix in Figure~\ref{figure_1D_H2_3level_extended_sparsity_reordered} to a conventional sparse matrix solver, \textbf{we will not obtain an $\mathcal{O}(N)$ algorithm}. This is because there will be fill-in's, which is detrimental to the linear scaling. For instance, if we eliminate the set of rows and columns corresponding to the leaf level without any compression, the matrix we obtain has a complete fill-in as shown in Figure~\ref{figure_1D_H2_2level_fillin_extended_sparsity_reordered}. For the linear scaling, what we need is that the matrix pattern after eliminating the rows and columns of leaf, should look like Figure~\ref{figure_1D_H2_2level_extended_sparsity_reordered}. The fill-ins are detrimental to the scaling of the algorithm. In our algorithm, these fill-ins are compressed since these correspond to interaction between well-separated clusters. Let us see how this is done by viewing the appropriate graph of this extended sparse-matrix.
\end{remark}

\begin{figure}[!htbp]
\includegraphics[scale=0.2]{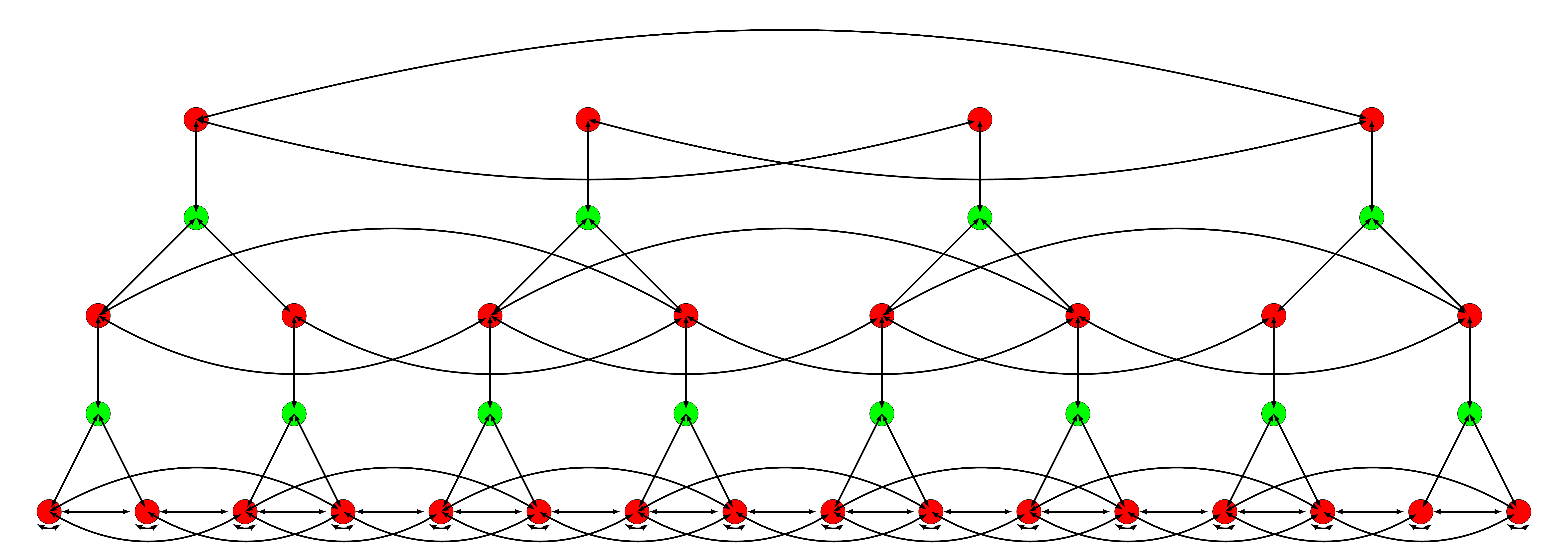}
\caption{The graph of the $3$ level extended sparse matrix in Figure~\ref{figure_1D_H2_3level_extended_sparsity_reordered}. The leaf consists of the bottom $16$ red nodes and the $8$ green nodes connected to them.}
\label{figure_extended_sparse_matrix_graph}
\end{figure}

\begin{figure}[!htbp]
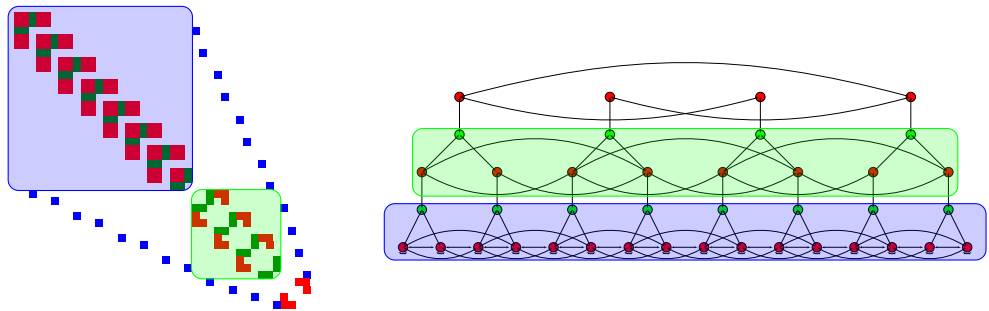

\tikz{
	\node at (0,0) {\includegraphics[scale=0.25]{./images/1D_H2_3level.pdf}};
	\node at (7,0) {\includegraphics[scale=0.125]{./images/elimination_matrix/elimination_matrix.pdf}};
	\fill [draw=blue, fill=blue, fill opacity=0.2, rounded corners] (-2,-0.45) rectangle (0.45,2);
	\fill [draw=blue, fill=blue, fill opacity=0.2, rounded corners] (3,-0.625) rectangle (11,-1.375);
	\fill [draw=green, fill=green, fill opacity=0.2, rounded corners] (0.4375,-0.4375) rectangle (1.625,-1.625);
	\fill [draw=green, fill=green, fill opacity=0.2, rounded corners] (3.375,-0.525) rectangle (10.625,0.375);
}
\caption{The duality between corresponding blocks of the matrix and the nodes \& edges of the graph at the leaf level are highlighted.}
\label{figure_graph_matrix_duality}
\end{figure}

\begin{figure}[!htbp]
\includegraphics[scale=0.2]{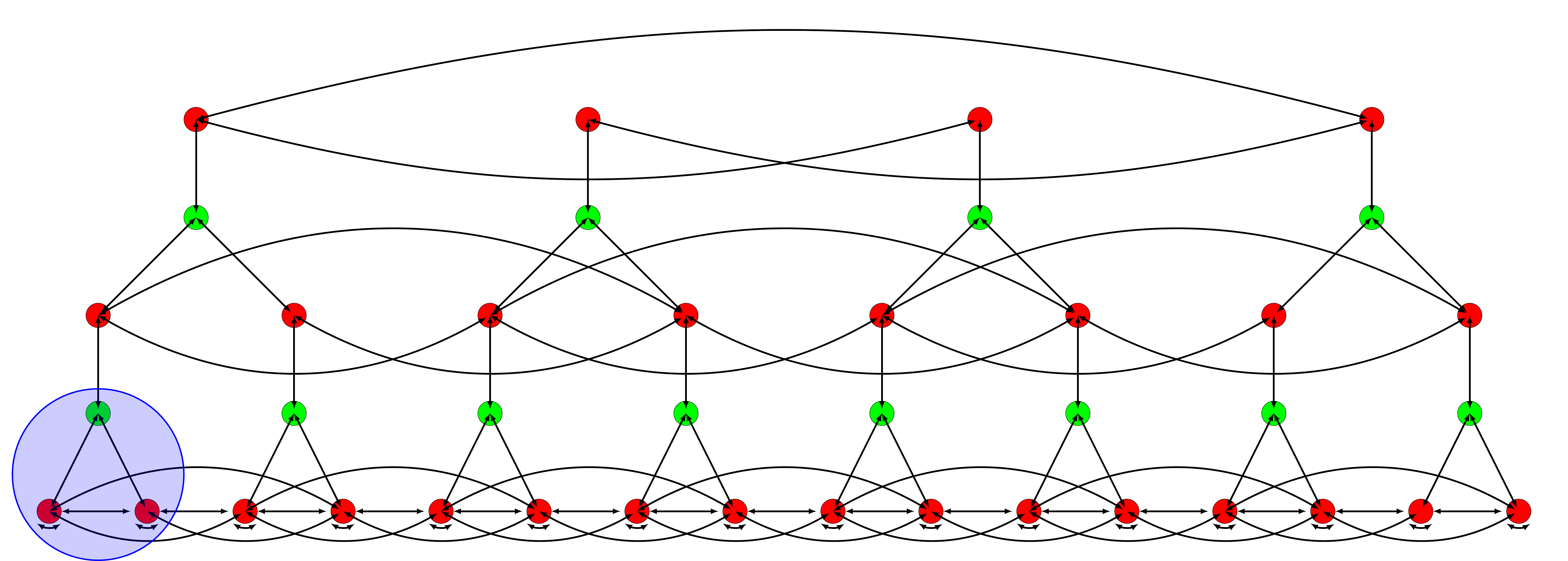}
\caption{Preparing to eliminate the particles and the corresponding local coefficients at the leaf level for the first set of clusters at the leaf level.}
\label{figure_elimination_matrix_first_step}
\end{figure}

\begin{figure}[!htbp]
\includegraphics[scale=0.2]{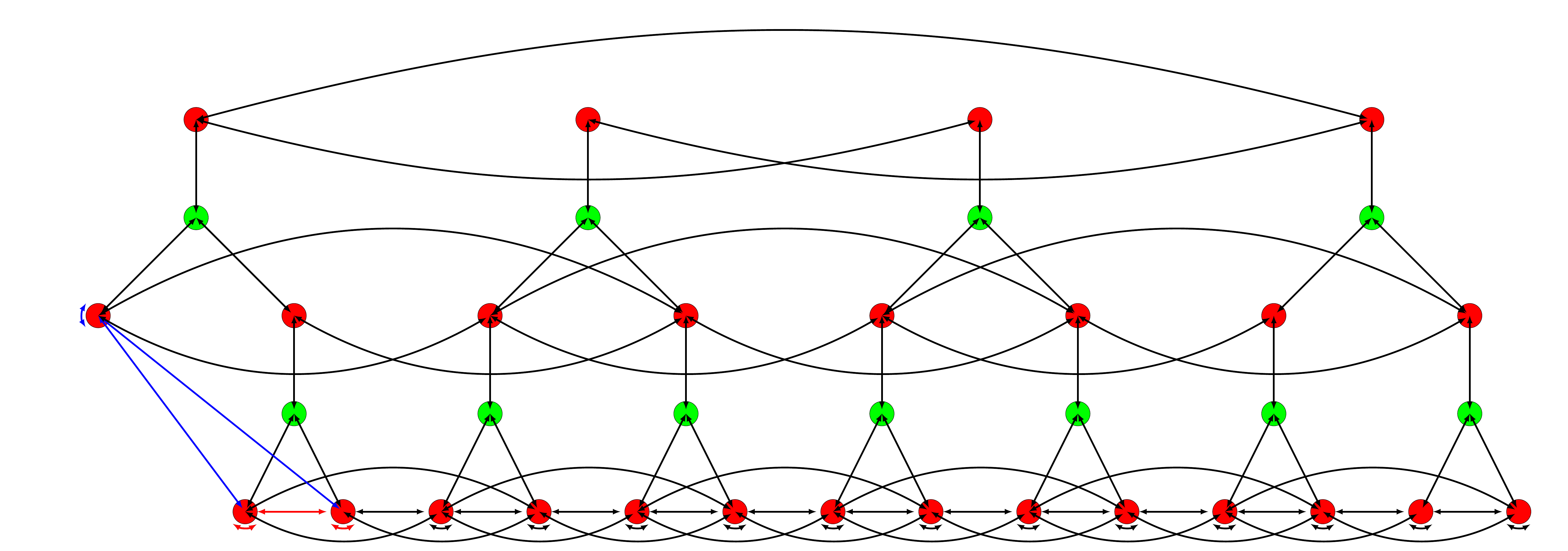}
\caption{Eliminating the first set of clusters in the previous step results in a fill-in (the edges corresponding to the fill-in are shown in blue color) between multipoles of the cluster with itself and with the particles of the neighbor. This also updates the interaction between the particles of the neighbors (the updated edges are shown in red color).}
\label{figure_elimination_matrix_after_first_step}
\end{figure}

\begin{figure}[!htbp]
\includegraphics[scale=0.2]{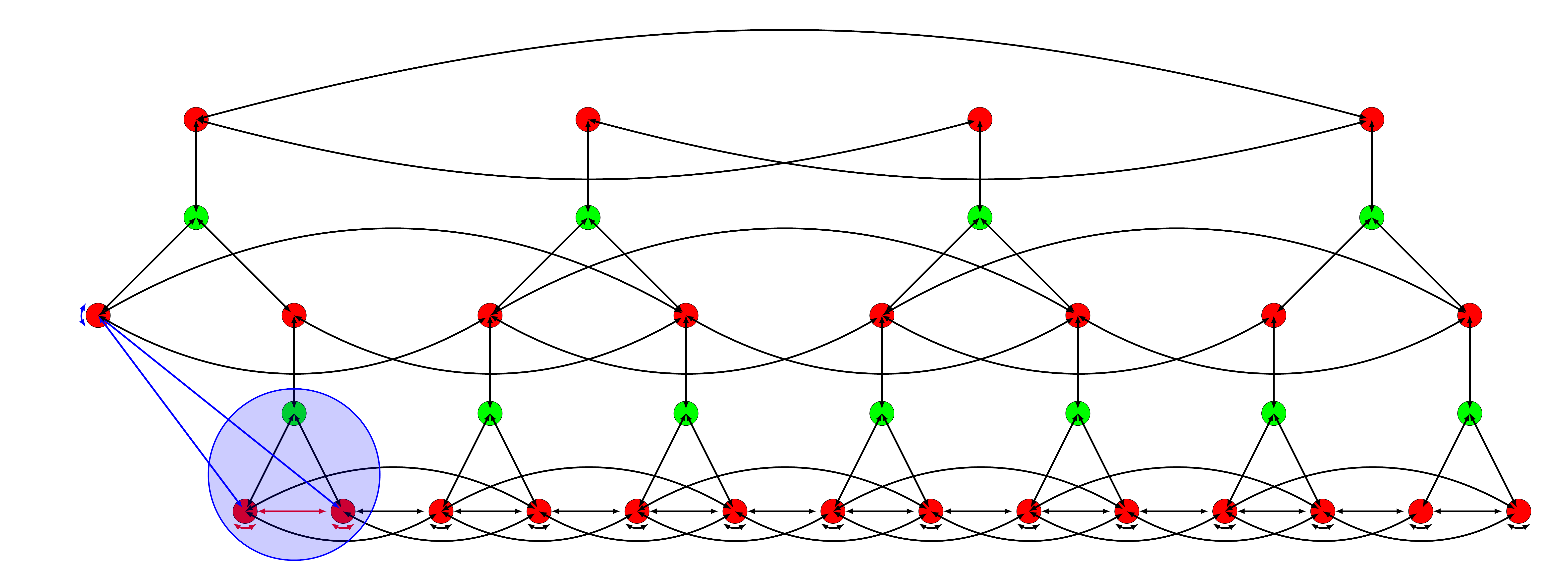}
\caption{Preparing to eliminate the next set of particles and the local coefficients of the second set of clusters at the leaf level.}
\label{figure_elimination_matrix_second_step}
\end{figure}

\begin{figure}[!htbp]
\includegraphics[scale=0.2]{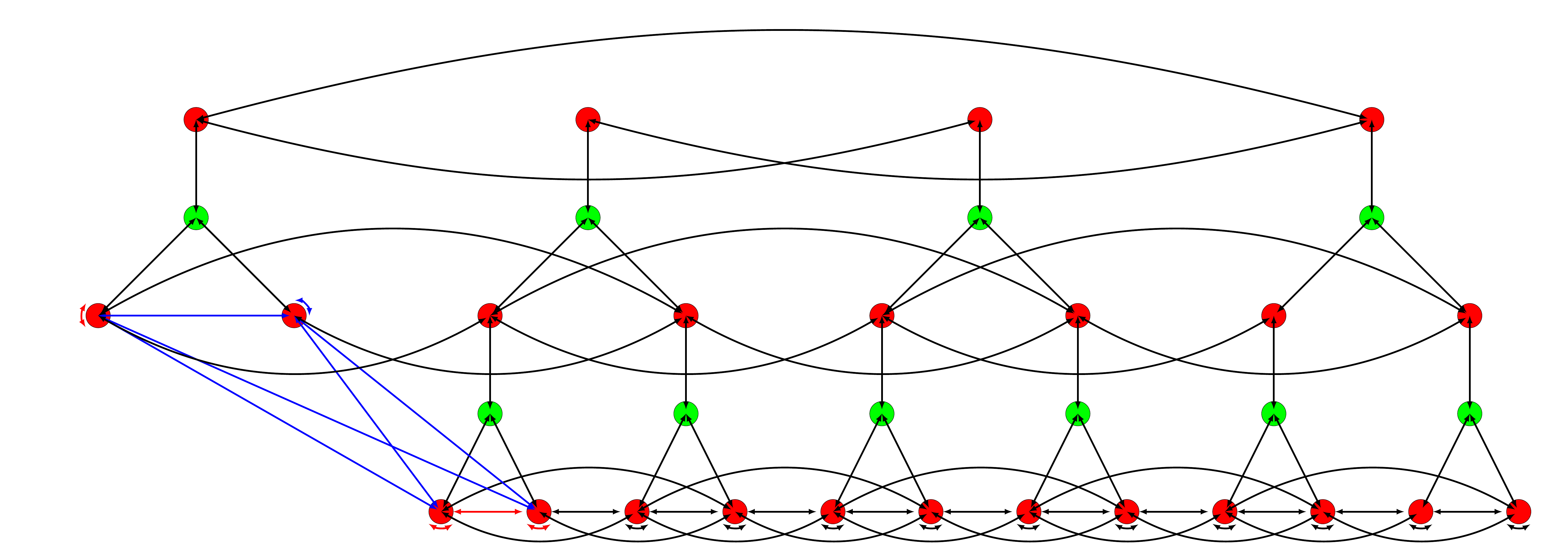}
\caption{Eliminating the second set of clusters result in fill-in between the multipoles of the first, multipoles of the second cluster and the particles of the third cluster, in addition to a fill-in in the multipoles of the second cluster with itself. Further, this also updates the interaction of the first set of multipoles with itself and the interaction of the particles of the third cluster with itself.}
\label{figure_elimination_matrix_after_second_step}
\end{figure}

\begin{figure}[!htbp]
\includegraphics[scale=0.2]{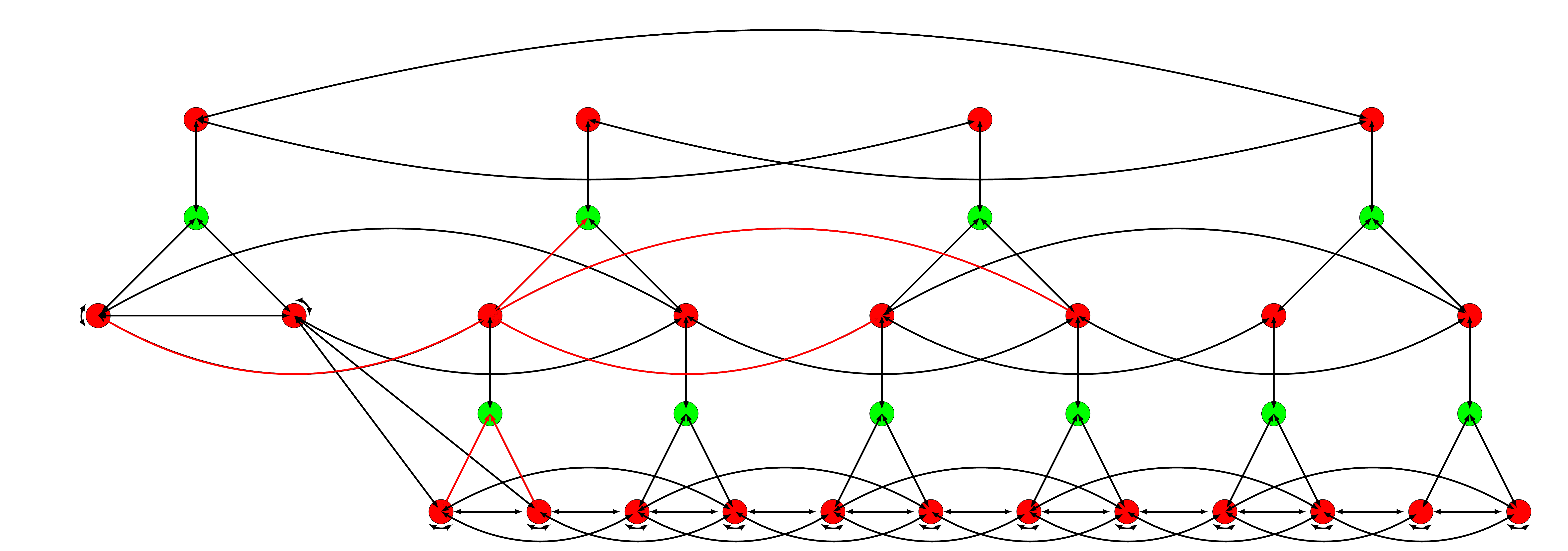}
\caption{This step is the key step since this reduces the fill-ins as we proceed with the elimination. As we saw in the previous step (Figure~\ref{figure_elimination_matrix_after_second_step}), the fill-in appears between the multipoles of the first set of clusters and the particles of the third set of clusters. This fill-in corresponds to a well-separated interaction. This fill-in is compressed along with the L2L/M2M basis of the third cluster, i.e., the interaction between the multipoles of the first cluster and the particles of the third cluster is redirected through the multipoles and local of the third cluster. This results in updating the following operators of the second cluster: (a) the M2L operator capturing the interaction with its interaction list, (b) the P2M and L2P operators, (c) the P2M and L2P operator of its parent.}
\label{figure_elimination_matrix_after_second_step_compression}
\end{figure}

\begin{figure}
\includegraphics[scale=0.2]{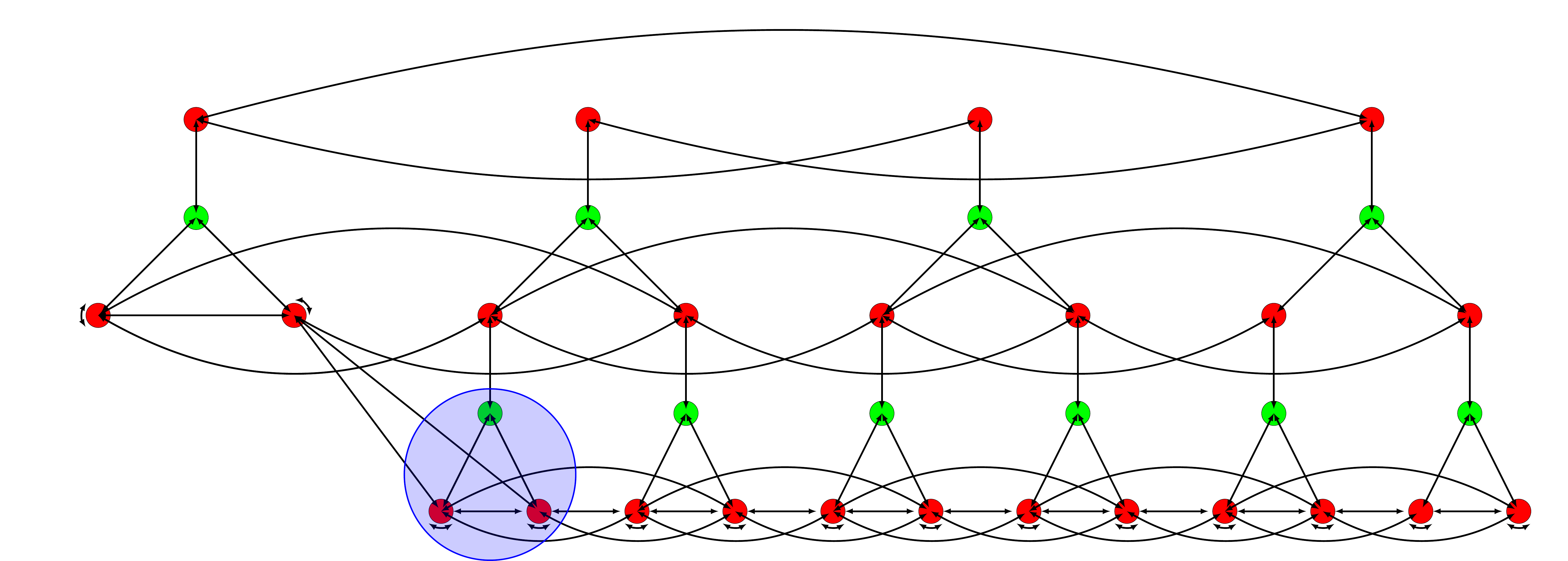}
\caption{Preparing to eliminate the third set of particles and their local coefficients at the leaf level.}
\label{figure_elimination_matrix_third_step}
\end{figure}

\begin{figure}
\includegraphics[scale=0.2]{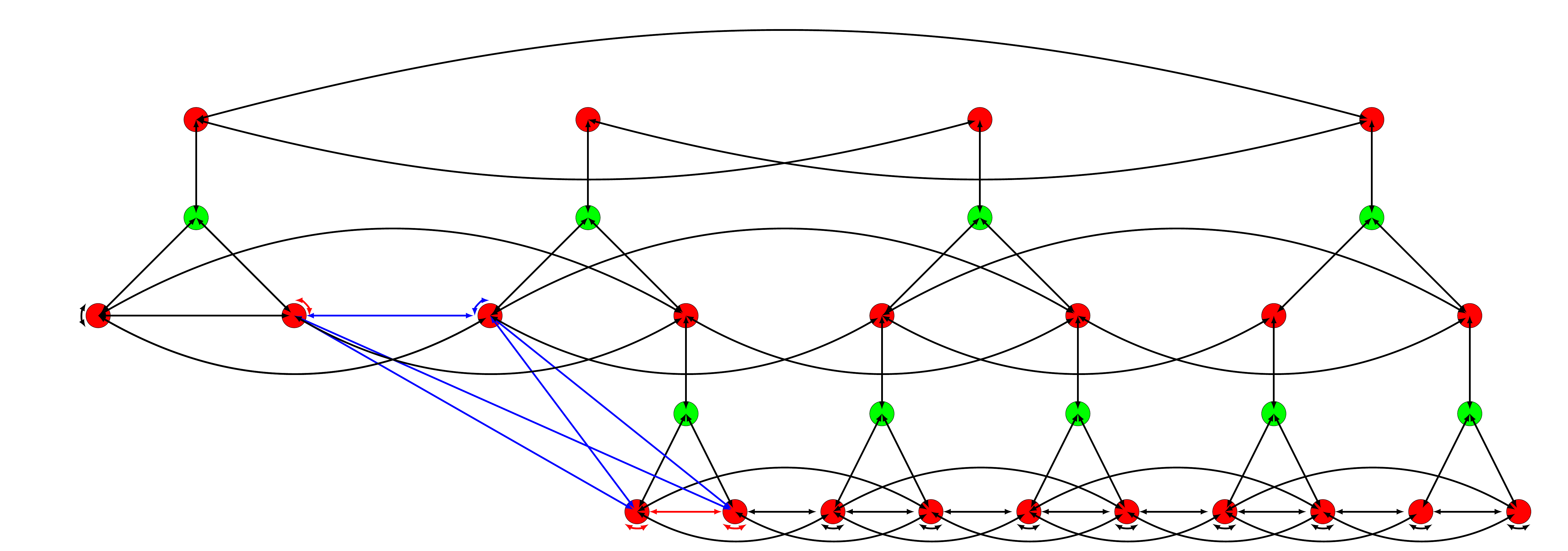}
\caption{Eliminating the third set of clusters result in fill-in between the multipoles of the second, multipoles of the third cluster and the particles of the fourth cluster, in addition to a fill-in in the multipoles of the third cluster with itself. Further, this also updates the interaction of the second set of multipoles with itself and the interaction of the particles of the fourth cluster with itself.}
\label{figure_elimination_matrix_after_third_step}
\end{figure}

\begin{figure}[!htbp]
\includegraphics[scale=0.2]{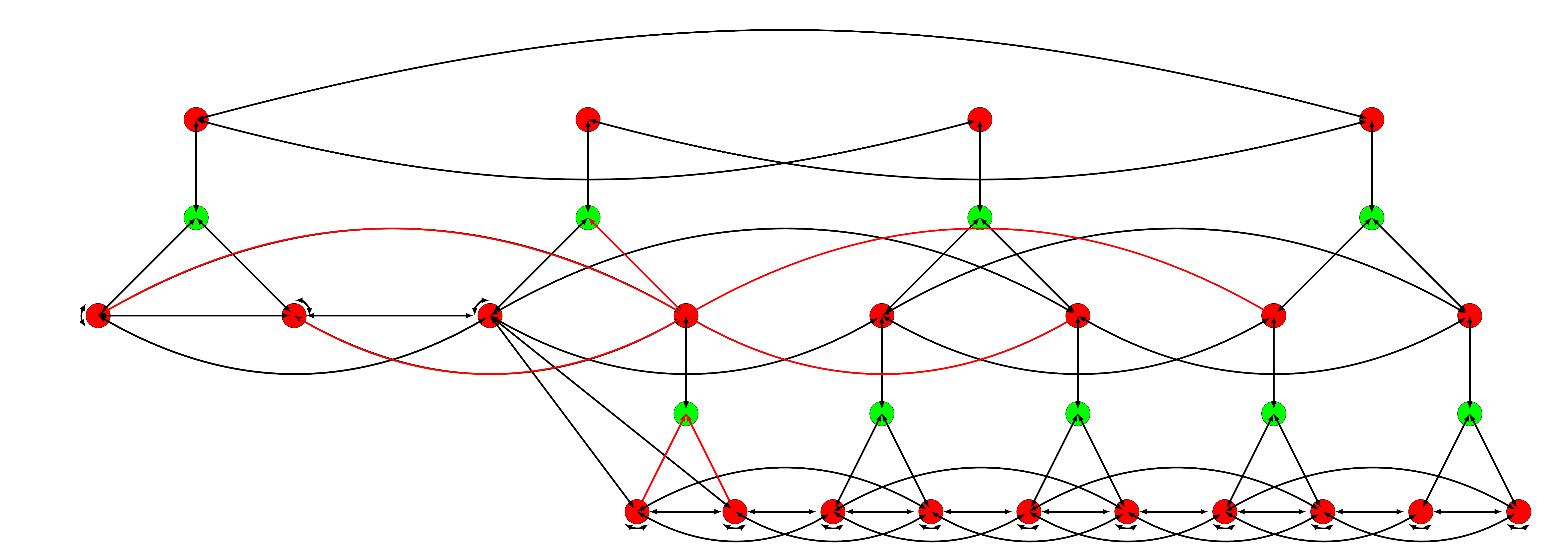}
\caption{This step reduces the fill-ins as we proceed with the elimination. As we saw in the previous step (Figure~\ref{figure_elimination_matrix_after_third_step}), the fill-in appears between the multipoles of the second set of clusters and the particles of the fourth set of clusters. This fill-in corresponds to a well-separated interaction. This fill-in is compressed along with the L2L/M2M basis of the fourth cluster, i.e., the interaction between the multipoles of the second cluster and the particles of the fourth cluster is redirected through the multipoles and local of the fourth cluster. This results in updating the following operators of the fourth cluster: (a) the M2L operator capturing the interaction with its interaction list, (b) the P2M and L2P operators, (c) the P2M and L2P operator of its parent.}
\label{figure_elimination_matrix_after_third_step_compression}
\end{figure}

Figures~\ref{figure_extended_sparse_matrix_graph}\textendash~\ref{figure_elimination_matrix_after_third_step_compression} explain the first few steps of the algorithm in $1$D.
\begin{remark}
In $2$D and $3$D, eliminating a cluster will also result in P2P fill-in. If this P2P fill-in corresponds to an interaction between well-separated clusters, this needs to be compressed as well. On an interval in $1$D, eliminating under the natural ordering results in no P2P fill-in. However, if we were to eliminate in order other than the natural ordering, there will be a P2P fill-in, which needs to be compressed.
\label{remark_1Ddiff}
\end{remark}

\begin{figure}[!htbp]
\includegraphics[scale=0.2]{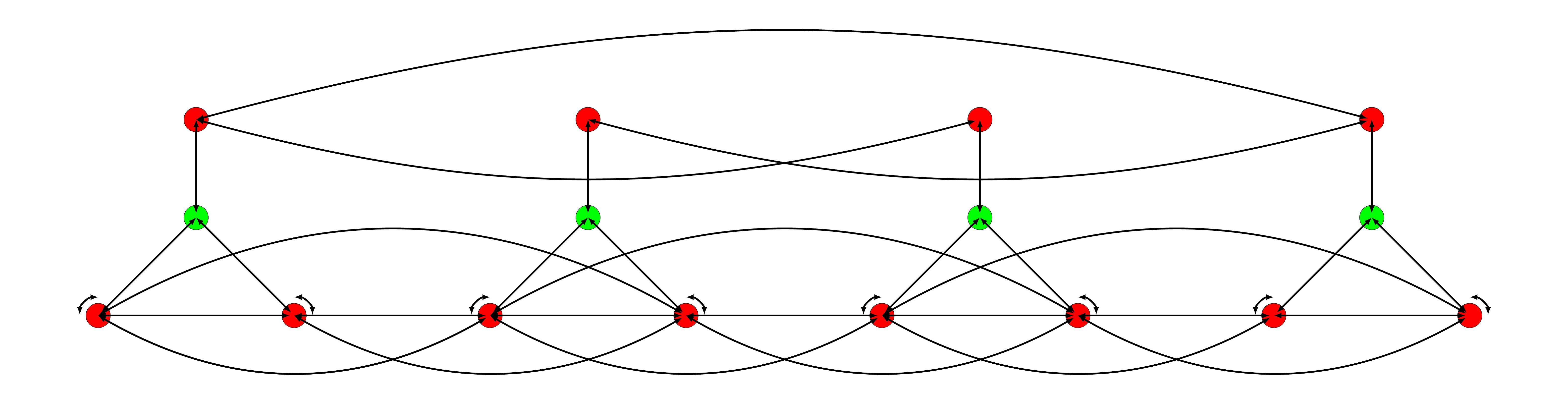}
\caption{Repeat this for all clusters at the leaf level (in our case the third level). Now repeat the process for the all the clusters at the second level by interpreting the multipoles of the leaf level as particles at this level. And proceed level by level. Once the elimination phase is done, we need to do the back substitution phase, which proceeds exactly in the reverse order.}
\label{figure_elimination_matrix_after_eliminating_leaf_levels}
\end{figure} 

\FloatBarrier

\subsection{General algorithm}
\label{subsection_Algorithm}
In general, in any dimension, set up the FMM tree, i.e., sub-divide the domain using a $2^d$ tree in $d$ dimensions. We will work with a non-adaptive/uniform tree for pedagogical reasons. At each level of the tree, we have unknown particles, multipoles and locals. We first introduce some notations in Table~\ref{table_cluster} to denote different clusters.
\begin{table}[!htbp]
\caption{Different lists associated with cluster `$i$' at level `$k$'}
\rowcolors{1}{}{gray!40}
\begin{tabular}{|c|l|}
\hline
$\mathcal{C}_i^{(k)}$ & Cluster $i$ at level $k$.\\
\hline
$\mathcal{N}_i^{(k)}$ & Neighbors of $\mathcal{C}_i^{(k)}$ including self.\\
\hline
{$\mathcal{I}_i^{(k)}$} & Interaction list, i.e., well-separated clusters\\
& that are children of parents neighbors.\\
\hline
\end{tabular}
\label{table_cluster}
\end{table} 

Each cluster $\mathcal{C}_i^{(k)}$ has the variables shown in Table~\ref{table_variables}.

\begin{table}[!htbp]
\caption{Known and unknown variables associated with cluster `$i$' at level `$k$'}
\rowcolors{1}{}{gray!40}
\begin{tabular}{|c|l|}
\hline
$x_i^{(k)}$ & Unknown charges on particles\\
\hline
$y_i^{(k)}$ & Unknown multipole coefficients\\
\hline
$z_i^{(k)}$ & Unknown local coefficients\\
\hline
$r_i^{(k)}$ & Known right hand side or the potential\\
\hline
$E_i^{(k)}$ & $1$ if the cluster is eliminated; $0$ otherwise\\
\hline
\end{tabular}
\label{table_variables}
\end{table}

Note that by particles at a non-leaf level, we mean the multipoles of its children, i.e.,
\begin{align}
x_i^{(k)} =
\begin{bmatrix}
y_{i_1}^{(k+1)} & y_{i_2}^{(k+1)} & \cdots & y_{i_{2^d}}^{(k+1)}
\end{bmatrix}^T
\label{equation_particles_multipoles}
\end{align}
where $i_p$ is a child of cluster $i$. As for the right hand side, the right hand side at each of the non-leaf level is set as zero, while the right hand side for each cluster at the leaf level is the input right hand side, i.e., the potential at these points. Refer Equation~\eqref{equation_1D_H2_2level_extended_sparsity_reordered} as to why this should be the case. The local coefficients $z_i^{(k)}$ is the potential on the multipoles due to well-separated clusters. Note that $z_i^{(1)} = 0$, since there is no well-separated cluster at the first level in the tree.

Table~\ref{table_FMM_operators} presents the operators needed for the FMM. In case of the FMM, $P2P_{ij}^{(k)}$ is non-zero only for $\mathcal{C}_i^{(k)} \in \mathcal{N}_j^{(k)}$ and $M2L_{ij}^{(k)}$ is non-zero for $\mathcal{C}_i^{(k)} \in \mathcal{I}_j^{(k)}$.
\begin{table}[!htbp]
\caption{Different interaction operators associated with cluster `$i$' at level `$k$'}
\rowcolors{1}{}{gray!40}
\begin{tabular}{|c|l|}
\hline
$P2P_{ij}^{(k)}$ & Potential on particles in $\mathcal{C}_i^{(k)}$ due to charges in $\mathcal{C}_j^{(k)}$.\\
\hline
$P2M_{ii}^{(k)}$ & Lumping the charges in $\mathcal{C}_i^{(k)}$ to its multipoles.\\
\hline
$L2P_{ii}^{(k)}$ & Interpolating the potential from locals in $\mathcal{C}_i^{(k)}$ to its particles.\\
\hline
$M2L_{ij}^{(k)}$ & Potential on locals of $\mathcal{C}_i^{(k)}$ due to the multipoles of $\mathcal{C}_j^{(k)}$.\\
\hline
\end{tabular}
\label{table_FMM_operators}
\end{table} 

For the IFMM, there will be a fill-in in $P2P_{ij}^{(k)}$, where $\mathcal{C}_i^{(k)} \in \mathcal{I}_j^{(k)}$, though this will be compressed on the fly. Apart from this, we need two additional operators for the IFMM as shown in Table~\ref{table_IFMM_operators}.

\begin{table}[!htbp]
\caption{Different interaction operators associated with cluster `$i$' at level `$k$'}
\rowcolors{1}{}{gray!40}
\begin{tabular}{|c|l|}
\hline
$P2L_{ij}^{(k)}$ & Local potentials in $\mathcal{C}_i^{(k)}$ due to charges in $\mathcal{C}_j^{(k)}$.\\
\hline
$M2P_{ij}^{(k)}$ & Potential on particles in $\mathcal{C}_i^{(k)}$ due to multipoles in $\mathcal{C}_j^{(k)}$.\\
\hline
\end{tabular}
\label{table_IFMM_operators}
\end{table} 

The $P2L_{ij}^{(k)}$ and $M2P_{ij}^{(k)}$ operators are needed only if $\mathcal{C}_i^{(k)} \in \mathcal{N}_j^{(k)} \cup \mathcal{I}_j^{(k)}$ and will be zero to begin with. However, while performing elimination these operators will get populated. The $P2L_{ij}^{(k)}$ and $M2P_{ij}^{(k)}$ operators for $\mathcal{C}_i^{(k)} \in \mathcal{I}_j^{(k)}$ will be compressed on the fly as well, i.e., for well-separated clusters
\begin{itemize}
\item
$P2L_{ij}^{(k)}$ will be expressed using the $M2L_{ij}^{(k)}$ and $P2M_{jj}^{(k)}$ by beefing up the ``multipoles" of $\mathcal{C}_j^{(k)}$.
\item
$M2P_{ij}^{(k)}$ will be expressed using the $L2P_{ii}^{(k)}$ and $M2L_{ij}^{(k)}$ by beefing up the ``locals" of $\mathcal{C}_i^{(k)}$.
\end{itemize}

\begin{remark}
$P2P$ operator at non-leaf levels is initially zero. As we proceed with the elimination, the P2P operator will be defined using the $M2L$ operator of its children, i.e.,
$$P2P_{ij}^{(k)} =
\begin{bmatrix}
M2L_{i_1j_1}^{(k+1)} & M2L_{i_1j_2}^{(k+1)} & M2L_{i_1j_3}^{(k+1)} & \cdots & M2L_{i_1j_{2^d}}^{(k+1)}\\
M2L_{i_2j_1}^{(k+1)} & M2L_{i_2j_2}^{(k+1)} & M2L_{i_2j_3}^{(k+1)} & \cdots & M2L_{i_2j_{2^d}}^{(k+1)}\\
M2L_{i_3j_1}^{(k+1)} & M2L_{i_3j_2}^{(k+1)} & M2L_{i_3j_3}^{(k+1)} & \cdots & M2L_{i_3j_{2^d}}^{(k+1)}\\
\vdots & \vdots & \vdots & \ddots & \vdots\\
M2L_{i_{2^d}j_1}^{(k+1)} & M2L_{i_{2^d}j_2}^{(k+1)} & M2L_{i_{2^d}j_3}^{(k+1)} & \cdots & M2L_{i_{2^d}j_{2^d}}^{(k+1)}\\
\end{bmatrix}
$$
where $i_p$ is a child of cluster $i$ and $j_q$ is a child of cluster $j$.
\end{remark}

\begin{remark}
Note that the M2L between two clusters at level $k+1$, i.e., the interaction between multipoles and locals of these two clusters, forms part of a neighbor interaction \textemdash a particle particle interaction \textemdash at level $k$.
\end{remark}
\begin{remark}
At the beginning of the algorithm, i.e., once the FMM data structure has been set up, the following are zero:
\begin{itemize}
\item
$P2L_{ij}^{(k)}, M2P_{ij}^{(k)}$;
\item
$M2L_{ij}^{(k)}$, where $j \in \mathcal{N}_i^{(k)}$;
\item
$P2P_{ij}^{(k)}$ and $r_i^{(k)}$ at the non-leaf levels;
\item
$E_j^{(k)} = 0$.
\end{itemize}

\end{remark}

\subsubsection{Overall idea of the algorithm}

Let $\kappa$ be the number of levels in the tree, i.e., the level $\kappa$ consists of leaves. Let $n_k$ denote the number of clusters at level $k$. For our purposes, we have $n_k = 2^k$. Below is a short snippet of the overall idea.
\begin{itemize}
\item
\textbf{Elimination phase/ Upward pass}
\begin{itemize}
\item
Eliminate cluster by cluster at the lowest level. Elimination order at level $k$: $\{e_1^{(k)}, e_2^{(k)}, \ldots, e_{n_{k}}^{(k)}\}$.
\item
Eliminating cluster, say $e_i^{(k)}$, results in fill-in among all its neighboring clusters, i.e., if $j_1, j_2 \in \mathcal{N}_i^{(k)}$, we then have a
\begin{itemize}
\item
$P2P_{j_1,j_2}^{(k)}$ and $P2P_{j_2,j_1}^{(k)}$ fill-in if $j_1,j_2$ have not been yet eliminated.
\item
$P2L_{j_1,j_2}^{(k)}$ and $M2P_{j_2,j_1}^{(k)}$ fill-in if $j_1$ has been eliminated and $j_2$ has not been yet eliminated.
\item
$M2L_{j_1,j_2}^{(k)}$ and $M2L_{j_2,j_1}^{(k)}$ fill-in if $j_1,j_2$ have been eliminated.
\end{itemize}
\item
Fill-ins between well-separated clusters, i.e., $P2P$, $P2L$, $M2P$ are compressed and directed through the appropriate $P2M$, $M2L$, $L2P$ operators, thereby eliminating the fill-ins.
\item
Repeat this at all levels marching up the tree, till we are left just with the particles at level $1$.
\item
Solve for these particles at level $1$, which are nothing but the multipoles at level $2$.
\end{itemize}
\item
\textbf{Back substitution phase/ Downward pass}
\begin{itemize}
\item
Back-substitute for each cluster starting from level $2$ marching downward in the tree.
\item
Back substitution order at level `$k$' is the reverse of the elimination order at level `$k$', i.e., $\{e_{n_k}^{(k)}, e_{n_k}^{(k-1)}, \ldots, e_{2}^{(k)}, e_{1}^{(k)}\}$
\end{itemize}
\end{itemize}

\subsubsection{Main algorithm in all dimensions}

The algorithm has an upward pass and a downward pass. The upward pass begins with clusters at the leaf level $\kappa$ and proceeds up till level $1$, while the downward pass begins with clusters at level $1$ and proceeds all the way till level $\kappa$. While performing the upward pass at level $k$, there are two main equations, P2P(i,k) and P2M(i,k), to be considered for cluster $i$ at this level.
\begin{itemize}
\item
\textbf{P2P equation}: Equation~\eqref{equation_p2p} will be denoted as P2P(i,k).
\begin{align}
\underbrace{P2P_{ii}^{(k)} x_i^{(k)}}_{\text{Direct self}} + \underbrace{L2P_{ii}^{(k)} z_{i}^{(k)}}_{\text{L2L self}} + \sum_{j \in \mathcal{N}_i^{(k)}} \left(\underbrace{P2P_{ij}^{(k)}x_j^{(k)} (1-E_j^{(k)})}_{\text{Direct neighbor}} + \underbrace{M2P_{ij}^{(k)} y_j^{(k)} E_j^{(k)}}_{\text{Multipoles of neighbors}}\right) = \begin{bmatrix} r_{i_1}^{(k+1)} \\ r_{i_2}^{(k+1)} \\ \cdots \\ r_{i_{2^d}}^{(k+1)}\end{bmatrix}
\label{equation_p2p}
\end{align}
\item
\textbf{P2M equation}: Equation~\eqref{equation_p2m} will be denoted as P2M(i,k).
\begin{align}
P2M_{ii}^{(k)} x_i^{(k)} - y_i^{(k)} = 0
\label{equation_p2m}
\end{align}
\end{itemize}

The elimination phase (upward pass) is presented in Algorithm~\ref{algorithm_elimination}, while Algorithm~\ref{algorithm_substitution} presents the back-substitution phase (downward pass). Note that the elimination phase involves more work (compressing the fill-ins) than the back-substitution phase. It is also important to note that, as with any direct solver, the ``factorization phase" (in our case most of the elimination phase), can be decoupled from the ``solve phase" (back-substitution phase) making it attractive for multiple right hand-sides.

\begin{algorithm}[!htbp]
\KwIn{The FMM tree}
\KwOut{Eliminates all the cluster}
\For{$k = \kappa, \kappa-1,\ldots,2$}{
	\For{$i = e_1^{(k)}, e_2^{(k)}, \ldots, e_{n_k}^{(k)}$}{
		Form $P2P$ operator for cluster $i$ with its neighbors from $M2L$ operators of its children\;
	}
	\For{$i = e_1^{(k)}, e_2^{(k)}, \ldots, e_{n_k}^{(k)}$}{
		Eliminate $x_i^{(k)}$ and $z_i^{(k)}$ using the equations P2P(i,k) and P2M(i,k)\;
		Set $E_i^{(k)} = 1$\;
		This results in the following fill-ins among neighbors of $i$\;
		Consider the pair $(p,q)$, where $\mathcal{C}_p^{(k)}, \mathcal{C}_q^{(k)} \in \mathcal{N}_i^{(k)}$\;
		\If{$p$ has been eliminated}{
			\If{$q$ has been eliminated}{
				Results in updating the $M2L_{pq}^{(k)}$ and $M2L_{qp}^{(k)}$ operator\;
			}
			\ElseIf{$q$ has not been eliminated}{
				Results in $P2L_{pq}^{(k)}$ and $M2P_{qp}^{(k)}$ fill-in\;
				\If{$p$ and $q$ are well-separated}{
					Compress $P2L_{pq}^{(k)}$, i.e., $P2L_{pq}^{(k)} = M2L_{pq}^{(k)} P2M_{qq}^{(k)}$\;
					Compress $M2P_{qp}^{(k)}$, i.e., $M2P_{qp}^{(k)} = L2P_{qq}^{(k)} M2L_{qp}^{(k)}$\;
					Update the following operators:
					\begin{itemize}
					\item
					$P2M_{qq}^{(k)}$, $L2P_{qq}^{(k)}$, $M2L_{lq}^{(k)}$, $M2L_{ql}^{(k)}$ operators, where $l \in \mathcal{N}_q^{(k)} \cup \mathcal{I}_q^{(k)}$\;
					\item
					$P2M_{q'q'}^{(k)}$, $L2P_{q'q'}^{(k)}$, where $q'$ is the parent of $q$\;
					\end{itemize}
				}	
			}
		}
		\ElseIf{$p$ has not been eliminated}{
			\If{$q$ has been eliminated}{
				Results in $P2L_{qp}^{(k)}$ and $M2P_{pq}^{(k)}$ fill-in\;
				\If{$p$ and $q$ are well-separated}{
					Compress $P2L_{qp}^{(k)}$, i.e., $P2L_{qp}^{(k)} = M2L_{qp}^{(k)} P2M_{pp}^{(k)}$\;
					Compress $M2P_{pq}^{(k)}$, i.e., $M2P_{pq}^{(k)} = L2P_{pp}^{(k)} M2L_{pq}^{(k)}$\;
					Update the following operators:
					\begin{itemize}
					\item
					$P2M_{pp}^{(k)}$, $L2P_{pp}^{(k)}$, $M2L_{lp}^{(k)}$, $M2L_{pl}^{(k)}$ operators, where $l \in \mathcal{N}_p^{(k)} \cup \mathcal{I}_p^{(k)}$\;
					\item
					$P2M_{p'p'}^{(k)}$, $L2P_{p'p'}^{(k)}$, where $p'$ is the parent of $p$\;
					\end{itemize}
				}
			}
			\ElseIf{$q$ has not been eliminated}{
				Results in $P2P_{pq}^{(k)}$ and $P2P_{qp}^{(k)}$ fill-in\;
				\If{$p$ and $q$ are well-separated}{
					Compress $P2P_{pq}^{(k)}$, i.e., $P2P_{pq}^{(k)} = L2P_{pp}^{(k)} M2L_{pq}^{(k)} P2M_{qq}^{(k)}$\;
					Compress $P2P_{qp}^{(k)}$, i.e., $P2P_{qp}^{(k)} = L2P_{qq}^{(k)} M2L_{qp}^{(k)} P2M_{pp}^{(k)}$\;
					Update the following operators:
					\begin{itemize}
					\item
					$L2P_{pp}^{(k)}$, $P2M_{pp}^{(k)}$, $L2P_{qq}^{(k)}$, $P2M_{qq}^{(k)}$\;
					\item
					$M2L_{lp}^{(k)}$, $M2L_{pl}^{(k)}$ operators, where $l \in \mathcal{N}_p^{(k)} \cup \mathcal{I}_p^{(k)}$\;
					\item
					$M2L_{lq}^{(k)}$, $M2L_{ql}^{(k)}$ operators, where $l \in \mathcal{N}_q^{(k)} \cup \mathcal{I}_q^{(k)}$\;
					\item
					$P2M_{p'p'}^{(k)}$,  $L2P_{p'p'}^{(k)}$, where $p'$ is the parent of $p$;
					\item
					$P2M_{q'q'}^{(k)}$,  $L2P_{q'q'}^{(k)}$, where $q'$ is the parent of $q$;
					\end{itemize}
				}	
			}
		}
	}
}
\caption{Elimination phase/ Upward pass}
\label{algorithm_elimination}
\end{algorithm}

\begin{algorithm}
\KwIn{Eliminated tree}
\KwOut{Back-Substituted tree}
\For{$k=1,2,\ldots,\kappa$}{
	\For{$i=e_{n_k}^{(k)}, e_{n_k-1}^{(k)}, \ldots, e_2^{(k)}, e_1^{(k)}$}{
		Obtain $x_i^{(k)}$ and $z_i^{(k)}$ using equation P2P(i,k) and P2M(i,k)\;
		Set $E_i^{(k)} = 0$\;
		Get $y_{i_1}^{(k+1)}, y_{i_2}^{(k+1)}, \ldots, y_{i_{2^d}}^{(k+1)}$ from $x_i^{(k)}$ using Equation~\eqref{equation_particles_multipoles}.
	}
}
\caption{Back-Substitution phase/ Downward pass}
\label{algorithm_substitution}
\end{algorithm}

\FloatBarrier

\section{Numerical benchmarks}

We present numerical benchmarks of the algorithm on a $2$D manifold for the following equation
\begin{align}
\sigma_i + \sum_{\overset{j=1}{j \neq i}}^{j =N}K(\vec{r}_i, \vec{r}_j) \sigma_j = f_i, \,\,\,\,\,\,\, \forall i \in \{1,2,\ldots,N\}
\label{main_Solving_Equation}
\end{align}
on three different singular kernels;

\begin{enumerate}[(i)]
\item
$\left( \dfrac{r(\log(r)-1)}{a(\log(a)-1)} \right) \chi_{r < a} + \left( \dfrac{\ln (r)}{\ln(a)} \right) \chi_{r \geq a}$
\item
$\left(\dfrac{r}a \right) \chi_{r<a} + \left( \dfrac{a}r \right) \chi_{r \geq a}$
\item
$\left( \dfrac{r^3(3\log(r)-1)}{a^3(3\log(a)-1)} \right) \chi_{r<a} + \left(\dfrac{r^2 \ln(r)}{a^2 \ln(a)} \right) \chi_{r \geq a}$
\end{enumerate}

We also provide comparison of the new solver with the HODLR fast direct solver discussed in~\cite{ambikasaran2013fast}, which is available here~\cite{ambikasaran2013HODLR} and also with the conventional full pivoted LU direct solver in Eigen. The algorithm was implemented in C++ and all the tests were run on $2.1$GHz with $256$GB memory. For each of the benchmark, we follow the conventions shown in Table~\ref{table_notations}.
\begin{table}[!htbp]
\caption{Notations used in the results}
\begin{tabular}{|c|l|}
\hline
$N$ & Number of unknowns\\
\hline
$r_m$ & Maximum rank of compressed sub-matrices\\
\hline
$t_a$ & Time taken to assemble the system\\
\hline
$t_f$ & Time taken to factor the system\\
\hline
$t_s$ & Time taken to apply the factorization, i.e., to solve the system\\
\hline
\multirow{3}{*}{Error} & A known vector $x_{\text{exact}}$ is taken and the right hand side $b$ is obtained.\\
& For this right hand side, the system is solved using the proposed algorithm\\
& and the relative error in $\ell_2$ norm is obtained, i.e., $\dfrac{\Vert x-x_{\text{exact}}\Vert_2}{\Vert x_{\text{exact}}\Vert_2}$\\
\hline
\end{tabular}
\label{table_notations}
\end{table}

The points $\vec{r}_i$ are distributed randomly in the square $[-1,1] \times [-1,1]$ such that the FMM tree is balanced. The nested low-rank decomposition of the well-separated clusters are obtained using Chebyshev interpolation by using $8$ Chebyshev nodes along one dimension, i.e., a total of $64$ Chebyshev nodes in $2$D, followed by SVD compression to further reduce the rank. The compression of the fill-ins is obtained using a tweaked version of the adaptive cross approximation algorithm~\cite{rjasanow2002adaptive,zhao2005adaptive}. The tolerance used to compress all the blocks is $10^{-14}$. The parameter $a$ was taken as $0.001$.

\begin{remark}
As seen from all the benchmarks, it is important to note that the rank of the compressed blocks for the IFMM remains independent of $N$, whereas for the HODLR solver the rank scales up (roughly like $\sqrt{N}$) with $N$.
\end{remark}

\begin{bm}
$K(r) = \left( \dfrac{r(\log(r)-1)}{a(\log(a)-1)} \right) \chi_{r < a} + \left( \dfrac{\ln (r)}{\ln(a)} \right) \chi_{r \geq a}$.
\end{bm}
\begin{table}[!htbp]
\caption{$K(r) = \left( \dfrac{r(\log(r)-1)}{a(\log(a)-1)} \right) \chi_{r < a} + \left( \dfrac{\ln (r)}{\ln(a)} \right) \chi_{r \geq a}$, where points are distributed in $[-1,1]^2$}
\resizebox{\hsize}{!}{
\begin{tabular}{|c|c|c|c|c|c|c|c|c|c|c|c|c|c|c|c|c|c|c|c|}
\hline
$N$ & \multicolumn{2}{c|}{$t_a$ in secs} & \multicolumn{2}{c|}{$t_f$ in secs} & \multicolumn{2}{c|}{$t_s$ in secs} & \multicolumn{2}{c|}{$r_m$} & \multicolumn{2}{c|}{Error}\\
\hline
(in thousands) & IFMM & HODLR & IFMM & HODLR & IFMM & HODLR & IFMM & HODLR & IFMM & HODLR\\
\hline
$1$ & $0.02$ & $0.07$ & $0.35$ & $0.09$ & $0.007$ & $0.0023$ & $29$ & $105$ & $10^{-12}$ & $10^{-11}$\\
\hline
$2$ & $0.02$ & $0.23$ & $0.71$ & $0.26$ & $0.009$ & $0.0055$ & $31$ & $137$ & $10^{-13}$ & $10^{-11}$\\
\hline
$5$ & $0.09$ & $1.24$ & $1.68$ & $1.21$ & $0.026$ & $0.0179$ & $30$ & $206$ & $10^{-12}$ & $10^{-9}$\\
\hline
$10$ & $0.12$ & $3.81$ & $3.21$ & $3.57$ & $0.06$ & $0.0467$ & $37$ & $257$ & $10^{-10}$ & $10^{-8}$\\
\hline
$20$ & $0.26$ & $12.15$ & $6.78$ & $11.31$ & $0.13$ & $0.1219$ & $36$ & $339$ & $10^{-11}$ & $10^{-8}$\\
\hline
$50$ & $0.71$ & $60.73$ & $16.79$ & $50.91$ & $0.32$ & $0.4237$ & $40$ & $485$ & - & -\\
\hline
$100$ & $1.78$ & $208.62$ & $34.93$ & $217.83$ & $0.72$ & $2.3386$ & $44$ & $671$ & - & -\\
\hline
$200$ & $3.43$ & - & $73.87$ & - & $1.54$ & - & $48$ & - & - & -\\ 
\hline
$500$ & $9.87$ & - & $189.63$ & - & $4.09$ & - & $47$ & - & - & -\\ 
\hline
$1000$ & $18.79$ & - & $412.48$ & - & $8.87$ & - & $51$ & - & - & -\\
\hline
\end{tabular}
}
\end{table}

\begin{figure}
\subfigure[Time taken for the solver]{
\includegraphics[scale=1.1625]{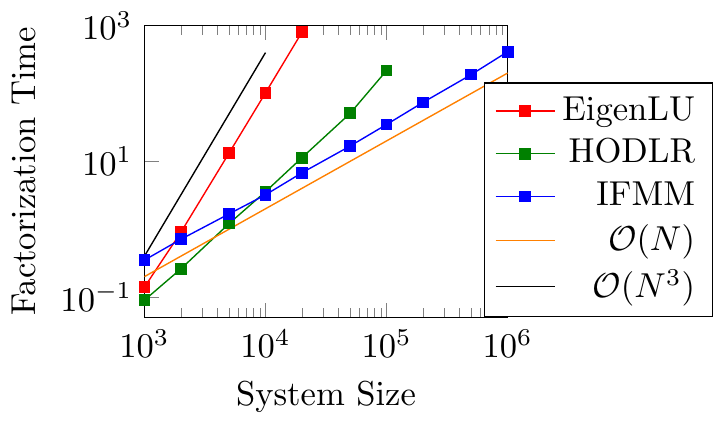}
}
\subfigure[Rank of the compressed blocks]{
\includegraphics[scale=1.1625]{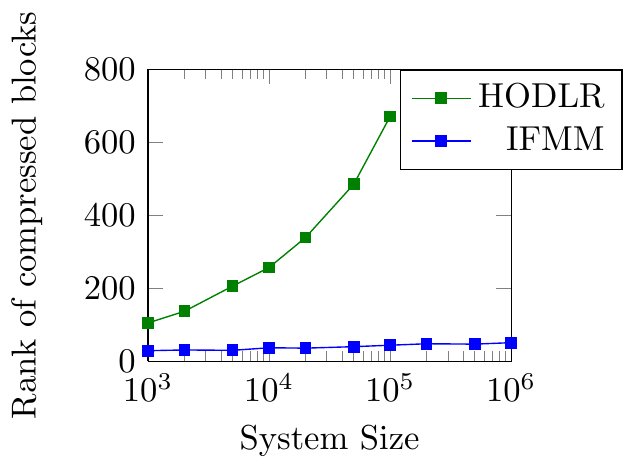}
}
\caption{$K(r) = \left( \dfrac{r(\log(r)-1)}{a(\log(a)-1)} \right) \chi_{r < a} + \left( \dfrac{\ln (r)}{\ln(a)} \right) \chi_{r \geq a}$, where points are distributed in $[-1,1]^2$}
\end{figure}
\FloatBarrier

\begin{bm}
$K(r) = \left( \dfrac{r^3(3\log(r)-1)}{a^3(3\log(a)-1)} \right) \chi_{r<a} + \left(\dfrac{r^2 \ln(r)}{a^2 \ln(a)} \right) \chi_{r \geq a}$.
\end{bm}
\begin{table}[!htbp]
\caption{$K(r) = \left( \dfrac{r^3(3\log(r)-1)}{a^3(3\log(a)-1)} \right) \chi_{r<a} + \left(\dfrac{r^2 \ln(r)}{a^2 \ln(a)} \right) \chi_{r \geq a}$, where points are distributed in $[-1,1]^2$}
\resizebox{\hsize}{!}{
\begin{tabular}{|c|c|c|c|c|c|c|c|c|c|c|c|c|c|c|c|c|c|c|c|}
\hline
$N$ & \multicolumn{2}{c|}{$t_a$ in secs} & \multicolumn{2}{c|}{$t_f$ in secs} & \multicolumn{2}{c|}{$t_s$ in secs} & \multicolumn{2}{c|}{$r_m$} & \multicolumn{2}{c|}{Error}\\
\hline
(in thousands) & IFMM & HODLR & IFMM & HODLR & IFMM & HODLR & IFMM & HODLR & IFMM & HODLR\\
\hline
$1$ & $0.01$ & $0.12$ & $0.11$ & $0.15$ & $0.003$ & $0.0029$ & $22$ & $129$ & $10^{-13}$ & $10^{-11}$\\
\hline
$2$ & $0.02$ & $0.36$ & $0.27$ & $0.42$ & $0.007$ & $0.0071$ & $21$ & $153$ & $10^{-12}$ & $10^{-10}$\\
\hline
$5$ & $0.08$ & $1.73$ & $0.64$ & $1.89$ & $0.017$ & $0.0244$ & $25$ & $222$ & $10^{-13}$ & $10^{-10}$\\
\hline
$10$ & $0.14$ & $4.87$ & $1.42$ & $5.47$ & $0.041$ & $0.0684$ & $24$ & $254$ & $10^{-11}$ & $10^{-9}$\\
\hline
$20$ & $0.23$ & $15.41$ & $2.76$ & $16.52$ & $0.091$ & $0.1578$ & $27$ & $314$ & $10^{-12}$ & $10^{-8}$\\
\hline
$50$ & $0.62$ & $61.21$ & $5.34$ & $68.57$ & $0.262$ & $0.9643$ & $26$ & $411$ & - & -\\
\hline
$100$ & $1.34$ & $166.03$ & $10.89$ & $231.77$ & $0.554$ & $3.5219$ & $28$ & $456$ & - & -\\
\hline
$200$ & $2.59$ & - & $26.47$ & - & $1.192$ & - & $31$ & - & - & -\\ 
\hline
$500$ & $6.43$ & - & $72.37$ & - & $2.873$ & - & $30$ & - & - & -\\ 
\hline
$1000$ & $12.12$ & - & $153.93$ & - & $6.121$ & - & $32$ & - & - & -\\
\hline
\end{tabular}
}
\end{table}

\begin{figure}
\subfigure[Time taken for the solver]{
\includegraphics[scale=1.1625]{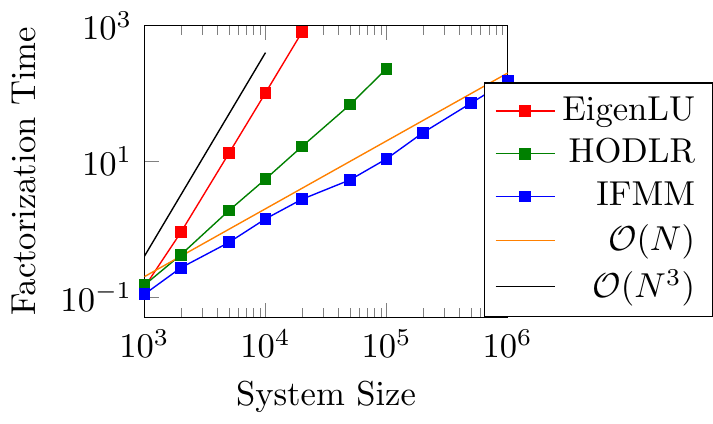}
}
\subfigure[Rank of the compressed blocks]{
\includegraphics[scale=1.1625]{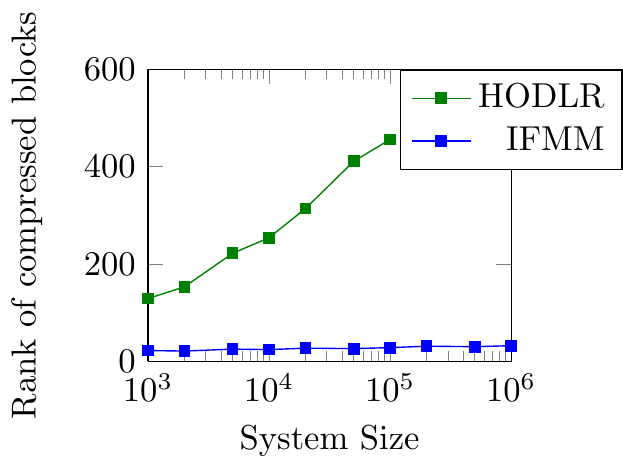}
}
\caption{$K(r) = \left( \dfrac{r^3(3\log(r)-1)}{a^3(3\log(a)-1)} \right) \chi_{r<a} + \left(\dfrac{r^2 \ln(r)}{a^2 \ln(a)} \right) \chi_{r \geq a}$, where points are distributed in $[-1,1]^2$}
\end{figure}
\FloatBarrier

\begin{bm}
$K(r) = \left(\dfrac{r}a \right) \chi_{r<a} + \left( \dfrac{a}r \right) \chi_{r \geq a}$.
\end{bm}

\begin{table}[!htbp]
\caption{$K(r) = \left(\dfrac{r}a \right) \chi_{r<a} + \left( \dfrac{a}r \right) \chi_{r \geq a}$, where points are distributed in $[-1,1]^2$}
\resizebox{\hsize}{!}{
\begin{tabular}{|c|c|c|c|c|c|c|c|c|c|c|c|c|c|c|c|c|c|c|c|}
\hline
$N$ & \multicolumn{2}{c|}{$t_a$ in secs} & \multicolumn{2}{c|}{$t_f$ in secs} & \multicolumn{2}{c|}{$t_s$ in secs} & \multicolumn{2}{c|}{$r_m$} & \multicolumn{2}{c|}{Error}\\
\hline
(in thousands) & IFMM & HODLR & IFMM & HODLR & IFMM & HODLR & IFMM & HODLR & IFMM & HODLR\\
\hline
$1$ & $0.02$ & $0.22$ & $0.42$ & $0.26$ & $0.005$ & $0.0042$ & $45$ & $141$ & $10^{-13}$ & $10^{-13}$\\
\hline
$2$ & $0.05$ & $0.79$ & $1.00$ & $0.62$ & $0.015$ & $0.0102$ & $43$ & $202$ & $10^{-12}$ & $10^{-12}$\\
\hline
$5$ & $0.11$ & $4.62$ & $2.37$ & $3.93$ & $0.035$ & $0.0490$ & $49$ & $252$ & $10^{-13}$ & $10^{-13}$\\
\hline
$10$ & $0.24$ & $18.93$ & $4.94$ & $16.26$ & $0.092$ & $0.1392$ & $51$ & $363$ & $10^{-11}$ & $10^{-11}$\\
\hline
$20$ & $0.51$ & $63.41$ & $9.92$ & $66.02$ & $0.192$ & $0.3478$ & $57$ & $501$ & $10^{-12}$ & $10^{-10}$\\
\hline
$50$ & $1.32$ & $332.3$ & $27.39$ & $415.84$ & $0.567$ & $2.0879$ & $54$ & $693$ & - & -\\
\hline
$100$ & $2.73$ & $1269.93$ & $61.34$ & $1712.79$ & $1.213$ & $7.9384$ & $59$ & $998$ & - & -\\
\hline
$200$ & $5.62$ & - & $142.38$ & - & $2.492$ & - & $65$ & - & - & -\\ 
\hline
$500$ & $14.52$ & - & $398.23$ & - & $5.781$ & - & $68$ & - & - & -\\ 
\hline
$1000$ & $28.42$ & - & $803.85$ & - & $12.436$ & - & $67$ & - & - & -\\
\hline
\end{tabular}
}
\end{table}

\begin{figure}
\subfigure[Time taken for the solver]{
\includegraphics[scale=1.1625]{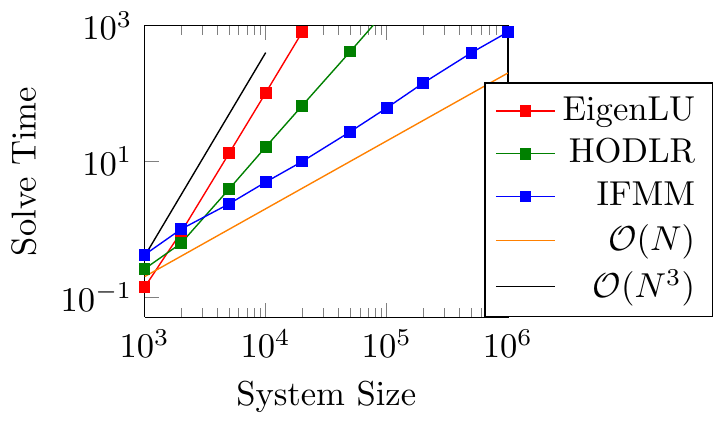}
}
\subfigure[Rank of the compressed blocks]{
\includegraphics[scale=1.1625]{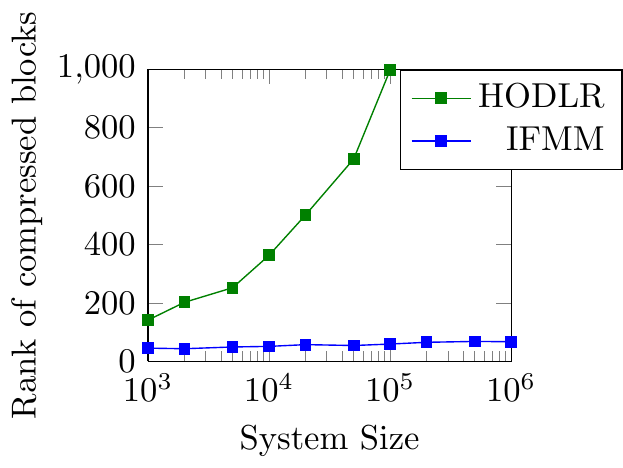}
}
\caption{$K(r) = \left(\dfrac{r}a \right) \chi_{r<a} + \left( \dfrac{a}r \right) \chi_{r \geq a}$, where points are distributed in $[-1,1]^2$}
\end{figure}
\FloatBarrier

\section{Conclusion}
\label{section_conclusion}
The article proposes the ``Inverse Fast Multipole Method". The IFMM is a fast direct solver that the solver works on the same data structure as the Fast Multipole Method and relies on compressing only the interactions corresponding to well-separated clusters. The highlight of the solver is that the computational cost scales linearly in the number of unknowns in all dimensions, provided the interactions corresponding to the well-separated clusters at all stages in the algorithm can be efficiently represented as a low-rank matrix. Numerical benchmarks presented validate the $\mathcal{O}(N)$ scaling of the algorithm for the kernels considered. The IFMM can be extended to integral equations, where the resulting linear system obtained after discretization can be solved at a computational complexity of $\mathcal{O}(N)$. It is also important to note that the IFMM can be applied to elliptic PDEs when discretized using local finite difference or fine element methods. In this case, the IFMM will operate on a sparse matrix, which is a special case of a hierarchical matrix with the rank corresponding to well-separated clusters being zero. The algorithm naturally extends itself to $\mathcal{H}^2$ matrices with strong admissibility criteria (weak low-rank structure).

\section{Acknowledgements}
Sivaram Ambikasaran would like to thank Leslie Greengard and Alex Barnett for helpful discussions in presenting the material. Sivaram Ambikasaran was supported by the Applied Mathematical Sciences Program of the U.S. Department of Energy under Contract DEFGO288ER25053 and by the Office of the Assistant Secretary of Defense for Research and Engineering and AFOSR under NSSEFF Program Award FA9550-10-1-0180. Part of this research was done at Stanford University, and was supported in part by the U.S.\ Army Research Laboratory, through the Army High Performance Computing Research Center, Cooperative Agreement W911NF-07-0027. This material is also based upon work supported by the Department of Energy National Nuclear Security Administration under Award Number {DE}-{NA}0002373-1.

\bibliographystyle{plainnat}

\begin{thebibliography}{66}
\providecommand{\natexlab}[1]{#1}
\providecommand{\url}[1]{\texttt{#1}}
\expandafter\ifx\csname urlstyle\endcsname\relax
  \providecommand{\doi}[1]{doi: #1}\else
  \providecommand{\doi}{doi: \begingroup \urlstyle{rm}\Url}\fi

\bibitem[{A}mbikasaran(2013)]{ambikasaran2013HODLR}
{S}ivaram {A}mbikasaran.
\newblock A fast direct solver for dense linear systems.
\newblock https://github.com/sivaramambikasaran/HODLR\_Solver, 2013.

\bibitem[Ambikasaran(2013)]{ambikasaran2013thesis}
Sivaram Ambikasaran.
\newblock \emph{Fast Algorithms for Dense Numerical Linear Algebra and
  Applications}.
\newblock PhD thesis, Stanford University, 2013.

\bibitem[Ambikasaran and Darve(2013)]{ambikasaran2013fast}
Sivaram Ambikasaran and Eric~F. Darve.
\newblock An $\mathcal{O}({N} \log {N})$ fast direct solver for partial
  hierarchically semi-separable matrices.
\newblock \emph{Journal of Scientific Computing}, 57\penalty0 (3):\penalty0
  477--501, 2013.

\bibitem[Ambikasaran and O'Neil(2014)]{ambikasaran2014fastsym}
Sivaram Ambikasaran and Michael O'Neil.
\newblock Fast symmetric factorization of hierarchical matrices with
  applications.
\newblock \emph{arXiv preprint arXiv:1405.0223}, 2014.

\bibitem[Ambikasaran et~al.(2013{\natexlab{a}})Ambikasaran, Li, Kitanidis, and
  Darve]{ambikasaran2013large}
Sivaram Ambikasaran, Judith~Y. Li, Peter~K. Kitanidis, and Eric~F. Darve.
\newblock Large-scale stochastic linear inversion using hierarchical matrices.
\newblock \emph{Computational Geosciences}, 17\penalty0 (6):\penalty0 913--927,
  2013{\natexlab{a}}.

\bibitem[Ambikasaran et~al.(2013{\natexlab{b}})Ambikasaran, Saibaba, Darve, and
  Kitanidis]{ambikasaran2013fastbayes}
Sivaram Ambikasaran, Arvind~K. Saibaba, Eric~F. Darve, and Peter~K. Kitanidis.
\newblock Fast algorithms for {B}ayesian inversion.
\newblock In \emph{Computational Challenges in the Geosciences}, pages
  101--142. Springer, 2013{\natexlab{b}}.

\bibitem[Ambikasaran et~al.(2014)Ambikasaran, Foreman-Mackey, Greengard, Hogg,
  and O'Neil]{ambikasaran2014fastdet}
Sivaram Ambikasaran, Daniel Foreman-Mackey, Leslie~F. Greengard, David~W. Hogg,
  and Michael O'Neil.
\newblock Fast direct methods for {G}aussian processes and the analysis of
  {NASA} {K}epler mission data.
\newblock \emph{arXiv preprint arXiv:1403.6015}, 2014.

\bibitem[Aminfar et~al.(2014)Aminfar, Ambikasaran, and Darve]{aminfar2014fast}
Amirhossein Aminfar, Sivaram Ambikasaran, and Eric~F. Darve.
\newblock A fast block low-rank dense solver with applications to
  finite-element matrices.
\newblock \emph{arXiv preprint arXiv:1403.5337}, 2014.

\bibitem[Arnoldi(1951)]{arnoldi1951principle}
Walter~E. Arnoldi.
\newblock The principle of minimized iterations in the solution of the matrix
  eigenvalue problem.
\newblock \emph{Quart. Appl. Math}, 9\penalty0 (1):\penalty0 17--29, 1951.

\bibitem[Barnes and Hut(1986)]{barnes1986hierarchical}
Josh Barnes and Piet Hut.
\newblock A hierarchical $\mathcal{O}({N} \log {N})$ force-calculation
  algorithm.
\newblock \emph{Nature}, 324\penalty0 (4):\penalty0 446--449, 1986.

\bibitem[Beatson and Greengard(1997)]{beatson1997short}
Rick~K. Beatson and Leslie~F. Greengard.
\newblock A short course on fast multipole methods.
\newblock \emph{Wavelets, multilevel methods and elliptic PDEs}, pages 1--37,
  1997.

\bibitem[Beatson and Newsam(1992)]{beatson1992fast}
Rick~K. Beatson and Garry~N. Newsam.
\newblock Fast evaluation of radial basis functions: {I}.
\newblock \emph{Computers \& Mathematics with Applications}, 24\penalty0
  (12):\penalty0 7--19, 1992.

\bibitem[Beatson et~al.(1999)Beatson, Cherrie, and Mouat]{beatson1999fast}
Rick~K. Beatson, Jon~B. Cherrie, and Cameron~T. Mouat.
\newblock Fast fitting of radial basis functions: Methods based on
  preconditioned {GMRES} iteration.
\newblock \emph{Advances in {C}omputational {M}athematics}, 11\penalty0
  (2):\penalty0 253--270, 1999.

\bibitem[Billings et~al.(2002)Billings, Beatson, and
  Newsam]{billings2002interpolation}
Stephen~D. Billings, Rick~K. Beatson, and Garry~N. Newsam.
\newblock Interpolation of geophysical data using continuous global surfaces.
\newblock \emph{Geophysics}, 67\penalty0 (6):\penalty0 1810, 2002.

\bibitem[B{\"o}rm(2006)]{borm2006matrix}
Steffen B{\"o}rm.
\newblock {$\mathcal{H}^2$}-matrix arithmetics in linear complexity.
\newblock \emph{Computing}, 77\penalty0 (1):\penalty0 1--28, 2006.

\bibitem[B\"{o}rm(2010)]{borm2010efficient}
Steffen B\"{o}rm.
\newblock Efficient numerical methods for non-local operators:
  $\mathcal{H}^2$-matrix compression, algorithms and analysis.
\newblock \emph{European Mathematical Society}, 14, 2010.

\bibitem[B{\"o}rm et~al.(2003)B{\"o}rm, Grasedyck, and
  Hackbusch]{borm2003introduction}
Steffen B{\"o}rm, Lars Grasedyck, and Wolfgang Hackbusch.
\newblock Introduction to hierarchical matrices with applications.
\newblock \emph{Engineering Analysis with Boundary Elements}, 27\penalty0
  (5):\penalty0 405--422, 2003.

\bibitem[Buhmann(2003)]{buhmann2003radial}
Martin~D. Buhmann.
\newblock \emph{Radial basis functions: theory and implementations}, volume~12.
\newblock Cambridge University Press, 2003.

\bibitem[Carr et~al.(2001)Carr, Beatson, Cherrie, Mitchell, Fright, McCallum,
  and Evans]{carr2001reconstruction}
J.C. Carr, Rick~K. Beatson, Jon~B. Cherrie, T.J. Mitchell, W.R. Fright, B.C.
  McCallum, and T.R. Evans.
\newblock Reconstruction and representation of {3D} objects with radial basis
  functions.
\newblock In \emph{Proceedings of the 28th annual conference on Computer
  graphics and interactive techniques}, pages 67--76. ACM, 2001.

\bibitem[Chandrasekaran et~al.(2002)Chandrasekaran, Dewilde, Gu, Pals, and
  van~der Veen]{chandrasekaran2002fast}
Shivkumar Chandrasekaran, Patrick Dewilde, Ming Gu, Timothy~P Pals, and
  AJ~van~der Veen.
\newblock \emph{Fast stable solver for sequentially semi-separable linear
  systems of equations}.
\newblock Springer, 2002.

\bibitem[Chandrasekaran et~al.(2006{\natexlab{a}})Chandrasekaran, Dewilde, Gu,
  Lyons, and Pals]{chandrasekaran2006fast1}
Shivkumar Chandrasekaran, Patrick Dewilde, Ming Gu, William Lyons, and
  Timothy~P Pals.
\newblock A fast solver for {HSS} representations via sparse matrices.
\newblock \emph{SIAM Journal on Matrix Analysis and Applications}, 29\penalty0
  (1):\penalty0 67--81, 2006{\natexlab{a}}.

\bibitem[Chandrasekaran et~al.(2006{\natexlab{b}})Chandrasekaran, Gu, and
  Pals]{chandrasekaran2006fast}
Shivkumar Chandrasekaran, Ming Gu, and Timothy~P Pals.
\newblock A fast {ULV} decomposition solver for hierarchically semi-separable
  representations.
\newblock \emph{SIAM Journal on Matrix Analysis and Applications}, 28\penalty0
  (3):\penalty0 603--622, 2006{\natexlab{b}}.

\bibitem[Cheng et~al.(1999)Cheng, Greengard, and Rokhlin]{cheng1999fast}
H.~Cheng, Leslie~F. Greengard, and Vladimir Rokhlin.
\newblock A fast adaptive multipole algorithm in three dimensions.
\newblock \emph{Journal of {C}omputational {P}hysics}, 155\penalty0
  (2):\penalty0 468--498, 1999.

\bibitem[Coifman et~al.(1993)Coifman, Rokhlin, and Wandzura]{coifman1993fast}
Ronald~R. Coifman, Vladimir Rokhlin, and S.~Wandzura.
\newblock The fast multipole method for the wave equation: A pedestrian
  prescription.
\newblock \emph{Antennas and Propagation Magazine, IEEE}, 35\penalty0
  (3):\penalty0 7--12, 1993.

\bibitem[Darve(2000{\natexlab{a}})]{darve2000fast}
Eric~F. Darve.
\newblock The fast multipole method: numerical implementation.
\newblock \emph{Journal of Computational Physics}, 160\penalty0 (1):\penalty0
  195--240, 2000{\natexlab{a}}.

\bibitem[Darve(2000{\natexlab{b}})]{darve2001fast}
Eric~F. Darve.
\newblock The fast multipole method {I}: Error analysis and asymptotic
  complexity.
\newblock \emph{SIAM Journal on Numerical Analysis}, 38\penalty0 (1):\penalty0
  98--128, 2000{\natexlab{b}}.

\bibitem[De~Boer et~al.(2007)De~Boer, Van~der Schoot, and Bijl]{de2007mesh}
A.~De~Boer, MS~Van~der Schoot, and H.~Bijl.
\newblock Mesh deformation based on radial basis function interpolation.
\newblock \emph{Computers \& {S}tructures}, 85\penalty0 (11-14):\penalty0
  784--795, 2007.

\bibitem[Fong and Darve(2009)]{fong2009black}
William Fong and Eric~F. Darve.
\newblock The black-box fast multipole method.
\newblock \emph{Journal of Computational Physics}, 228\penalty0 (23):\penalty0
  8712--8725, 2009.

\bibitem[Freund(1993)]{freund1993transpose}
Roland~W. Freund.
\newblock A transpose-free quasi-minimal residual algorithm for non-hermitian
  linear systems.
\newblock \emph{SIAM Journal on Scientific Computing}, 14:\penalty0 470, 1993.

\bibitem[Freund and Nachtigal(1991)]{freund1991qmr}
Roland~W. Freund and N\"{o}el~M. Nachtigal.
\newblock {QMR}: a quasi-minimal residual method for non-hermitian linear
  systems.
\newblock \emph{Numerische Mathematik}, 60\penalty0 (1):\penalty0 315--339,
  1991.

\bibitem[Grasedyck and Hackbusch(2003)]{grasedyck2003construction}
Lars Grasedyck and Wolfgang Hackbusch.
\newblock Construction and arithmetics of $\mathcal{H}$-matrices.
\newblock \emph{Computing}, 70\penalty0 (4):\penalty0 295--334, 2003.

\bibitem[Greengard(1988)]{greengard1988rapid}
Leslie~F. Greengard.
\newblock \emph{The rapid evaluation of potential fields in particle systems},
  volume 1987.
\newblock the MIT Press, 1988.

\bibitem[Greengard and Rokhlin(1987)]{greengard1987fast}
Leslie~F. Greengard and Vladimir Rokhlin.
\newblock A fast algorithm for particle simulations.
\newblock \emph{Journal of {C}omputational {P}hysics}, 73\penalty0
  (2):\penalty0 325--348, 1987.

\bibitem[Greengard and Rokhlin(1997)]{greengard1997new}
Leslie~F. Greengard and Vladimir Rokhlin.
\newblock A new version of the fast multipole method for the {L}aplace equation
  in three dimensions.
\newblock \emph{Acta {N}umerica}, 6\penalty0 (1):\penalty0 229--269, 1997.

\bibitem[Greengard et~al.(2009)Greengard, Gueyffier, Martinsson, and
  Rokhlin]{greengard2009fast}
Leslie~F. Greengard, Denis Gueyffier, Per-Gunnar Martinsson, and Vladimir
  Rokhlin.
\newblock Fast direct solvers for integral equations in complex
  three-dimensional domains.
\newblock \emph{Acta {N}umerica}, 18\penalty0 (1):\penalty0 243--275, 2009.

\bibitem[Gumerov and Duraiswami(2007)]{gumerov2007fast}
Nail~A. Gumerov and Ramani Duraiswami.
\newblock Fast radial basis function interpolation via preconditioned {K}rylov
  iteration.
\newblock \emph{SIAM Journal on Scientific Computing}, 29\penalty0
  (5):\penalty0 1876--1899, 2007.

\bibitem[Hackbusch(1999)]{hackbusch1999sparse}
Wolfgang Hackbusch.
\newblock A sparse matrix arithmetic based on $\mathcal{H}$-matrices. {P}art
  {I}: Introduction to $\mathcal{H}$-matrices.
\newblock \emph{Computing}, 62\penalty0 (2):\penalty0 89--108, 1999.

\bibitem[Hackbusch and B{\"o}rm(2002{\natexlab{a}})]{hackbusch2002data}
Wolfgang Hackbusch and Steffen B{\"o}rm.
\newblock Data-sparse approximation by adaptive $\mathcal{H}^2$-matrices.
\newblock \emph{Computing}, 69\penalty0 (1):\penalty0 1--35,
  2002{\natexlab{a}}.

\bibitem[Hackbusch and B{\"o}rm(2002{\natexlab{b}})]{hackbusch2002h2}
Wolfgang Hackbusch and Steffen B{\"o}rm.
\newblock $\mathcal{H}^2$-matrix approximation of integral operators by
  interpolation.
\newblock \emph{Applied Numerical Mathematics}, 43\penalty0 (1):\penalty0
  129--143, 2002{\natexlab{b}}.

\bibitem[Hackbusch and Khoromskij(2000{\natexlab{a}})]{hackbusch2000sparse}
Wolfgang Hackbusch and Boris~N Khoromskij.
\newblock A sparse $\mathcal{H}$-matrix arithmetic.
\newblock \emph{Computing}, 64\penalty0 (1):\penalty0 21--47,
  2000{\natexlab{a}}.

\bibitem[Hackbusch and Khoromskij(2000{\natexlab{b}})]{hackbusch2000sparse1}
Wolfgang Hackbusch and Boris~N Khoromskij.
\newblock A sparse $\mathcal{H}$-matrix arithmetic: general complexity
  estimates.
\newblock \emph{Journal of Computational and Applied Mathematics}, 125\penalty0
  (1):\penalty0 479--501, 2000{\natexlab{b}}.

\bibitem[Hackbusch and Nowak(1989)]{hackbusch1989fast}
Wolfgang Hackbusch and Z.P. Nowak.
\newblock On the fast matrix multiplication in the boundary element method by
  panel clustering.
\newblock \emph{Numerische Mathematik}, 54\penalty0 (4):\penalty0 463--491,
  1989.

\bibitem[Hackbusch et~al.(2000)Hackbusch, Khoromskij, and
  Sauter]{hackbusch2000h2}
Wolfgang Hackbusch, Boris Khoromskij, and Stefan~A Sauter.
\newblock On $\mathcal{H}^2$-matrices.
\newblock In Hans-Joachim Bungartz, Ronald~H.W. Hoppe, and Christoph Zenger,
  editors, \emph{Lectures on Applied Mathematics}, pages 9--29. Springer Berlin
  Heidelberg, 2000.
\newblock ISBN 978-3-642-64094-0.

\bibitem[Hackbusch et~al.(2001)Hackbusch, Grasedyck, and
  B{\"o}rm]{hackbusch2001introduction}
Wolfgang Hackbusch, Lars Grasedyck, and Steffen B{\"o}rm.
\newblock \emph{An introduction to hierarchical matrices}.
\newblock Max-Planck-Inst. f{\"u}r Mathematik in den Naturwiss., 2001.

\bibitem[Hager(1989)]{hager1989updating}
William~W. Hager.
\newblock Updating the inverse of a matrix.
\newblock \emph{SIAM review}, pages 221--239, 1989.

\bibitem[Hestenes and Stiefel(1952)]{hestenes1952methods}
Magnus~R. Hestenes and Eduard Stiefel.
\newblock Methods of conjugate gradients for solving linear systems.
\newblock \emph{Journal of Research of the National Bureau of Standards},
  49\penalty0 (6):\penalty0 409--436, 1952.

\bibitem[Ho and Greengard(2012)]{ho2012fast}
Kenneth~L. Ho and Leslie~F. Greengard.
\newblock A fast direct solver for structured linear systems by recursive
  skeletonization.
\newblock \emph{SIAM Journal on Scientific Computing}, 34\penalty0
  (5):\penalty0 2507--2532, 2012.

\bibitem[Kong et~al.(2011)Kong, Bremer, and Rokhlin]{kong2011adaptive}
Wai~Yip Kong, James Bremer, and Vladimir Rokhlin.
\newblock An adaptive fast direct solver for boundary integral equations in two
  dimensions.
\newblock \emph{{A}pplied and {C}omputational {H}armonic {A}nalysis},
  31\penalty0 (3):\penalty0 346--369, 2011.

\bibitem[Lai et~al.(2014)Lai, Ambikasaran, and Greengard]{lai2014fast}
Jun Lai, Sivaram Ambikasaran, and Leslie~F. Greengard.
\newblock A fast direct solver for high frequency scattering from a large
  cavity in two dimensions.
\newblock \emph{arXiv preprint arXiv:1404.3451}, 2014.

\bibitem[Lee et~al.(2013)Lee, Ambikasaran, Kitanidis, Illangasekare, and
  Smits]{lee2013hydrogeophysical}
Jonghyun Lee, Sivaram Ambikasaran, Peter~K. Kitanidis, Tissa~H. Illangasekare,
  and Kathleen~M. Smits.
\newblock Hydrogeophysical data assimilation using a fast {K}alman filter for
  managed aquifer recharge and recovery.
\newblock In \emph{AGU Fall Meeting Abstracts}, volume~1, page 1281, 2013.

\bibitem[Li et~al.(2014)Li, Ambikasaran, Darve, and Kitanidis]{li2014kalman}
Judith~Y. Li, Sivaram Ambikasaran, Eric~F. Darve, and Peter~K. Kitanidis.
\newblock A {K}alman filter powered by $\mathcal{H}^2$-matrices for
  quasi-continuous data assimilation problems.
\newblock \emph{Water Resources Research}, 2014.

\bibitem[Martinsson(2009)]{martinsson2009fast}
Per-Gunnar Martinsson.
\newblock A fast direct solver for a class of elliptic partial differential
  equations.
\newblock \emph{Journal of Scientific Computing}, 38\penalty0 (3):\penalty0
  316--330, 2009.

\bibitem[Martinsson and Rokhlin(2005)]{martinsson2005fast}
Per-Gunnar Martinsson and Vladimir Rokhlin.
\newblock A fast direct solver for boundary integral equations in two
  dimensions.
\newblock \emph{Journal of Computational Physics}, 205\penalty0 (1):\penalty0
  1--23, 2005.

\bibitem[Nishimura(2002)]{nishimura2002fast}
Naoshi Nishimura.
\newblock Fast multipole accelerated boundary integral equation methods.
\newblock \emph{Applied Mechanics Reviews}, 55\penalty0 (4):\penalty0 299--324,
  2002.

\bibitem[Paige and Saunders(1975)]{paige1975solution}
Christopher~C. Paige and Michael~A. Saunders.
\newblock Solution of sparse indefinite systems of linear equations.
\newblock \emph{SIAM Journal on Numerical Analysis}, 12\penalty0 (4):\penalty0
  617--629, 1975.

\bibitem[Pals(2004)]{pals2004multipole}
Timothy~P Pals.
\newblock \emph{Multipole for scattering computations: Spectral discretization,
  stabilization, fast solvers}.
\newblock PhD thesis, University of California Santa Barbara, 2004.

\bibitem[Rjasanow(2002)]{rjasanow2002adaptive}
Sergej Rjasanow.
\newblock Adaptive cross approximation of dense matrices.
\newblock \emph{IABEM 2002, International Association for Boundary Element
  Methods}, 2002.

\bibitem[Saad and Schultz(1986)]{saad1986gmres}
Yousef Saad and Martin~H. Schultz.
\newblock {GMRES}: A generalized minimal residual algorithm for solving
  nonsymmetric linear systems.
\newblock \emph{SIAM Journal on {S}cientific and {S}tatistical {C}omputing},
  7\penalty0 (3):\penalty0 856--869, 1986.

\bibitem[Saibaba et~al.(2012)Saibaba, Ambikasaran, Li, Kitanidis, and
  Darve]{saibaba2012application}
Arvind~K. Saibaba, Sivaram Ambikasaran, Judith~Y. Li, Peter~K. Kitanidis, and
  Eric~F. Darve.
\newblock Application of hierarchical matrices to linear inverse problems in
  geostatistics.
\newblock \emph{Oil and Gas Science and Technology-Revue de l'IFP-Institut
  Francais du Petrole}, 67\penalty0 (5):\penalty0 857, 2012.

\bibitem[Schaback(1995)]{schaback1995creating}
Robert Schaback.
\newblock Creating surfaces from scattered data using radial basis functions.
\newblock In \emph{Mathematical Methods for Curves and Surfaces}, pages
  477--496. University Press, 1995.

\bibitem[Van~der Vorst(1992)]{van1992bi}
Henk~A Van~der Vorst.
\newblock Bi-{CGSTAB}: A fast and smoothly converging variant of bi-{CG} for
  the solution of non-symmetric linear systems.
\newblock \emph{SIAM Journal on scientific and Statistical Computing},
  13\penalty0 (2):\penalty0 631--644, 1992.

\bibitem[Wang and Liu(2002)]{wang2002point}
J.~G. Wang and G.~R. Liu.
\newblock A point interpolation mesh-less method based on radial basis
  functions.
\newblock \emph{International Journal for Numerical Methods in Engineering},
  54\penalty0 (11):\penalty0 1623--1648, 2002.

\bibitem[Woodbury(1950)]{woodbury1950inverting}
Max~A Woodbury.
\newblock Inverting modified matrices.
\newblock Statistical Research Group, Memo. Rep. no. 42, Princeton University,
  1950.

\bibitem[Wu and Schaback(1993)]{wu1993local}
Zong-min Wu and Robert Schaback.
\newblock Local error estimates for radial basis function interpolation of
  scattered data.
\newblock \emph{{IMA} {J}ournal of {N}umerical {A}nalysis}, 13\penalty0
  (1):\penalty0 13--27, 1993.

\bibitem[Ying et~al.(2004)Ying, Biros, and Zorin]{ying2004kernel}
Lexing Ying, George Biros, and Denis Zorin.
\newblock A kernel-independent adaptive fast multipole algorithm in two and
  three dimensions.
\newblock \emph{Journal of Computational Physics}, 196\penalty0 (2):\penalty0
  591--626, 2004.

\bibitem[Zhao et~al.(2005)Zhao, Vouvakis, and Lee]{zhao2005adaptive}
Kezhong Zhao, Marinos~N Vouvakis, and J-F Lee.
\newblock The adaptive cross approximation algorithm for accelerated method of
  moments computations of {EMC} problems.
\newblock \emph{Electromagnetic Compatibility, IEEE Transactions on},
  47\penalty0 (4):\penalty0 763--773, 2005.

\end{thebibliography}

\end{document}